\documentclass{amsart}
\usepackage{amssymb,amsmath}

\usepackage[cmtip,all]{xy}
\usepackage{amscd}
\usepackage{graphics}

\usepackage[usenames]{color}
\definecolor{Red}{rgb}{1,0,0}

\newtheorem{theorem}{Theorem}[section]
\newtheorem{corollary}[theorem]{Corollary}
\newtheorem{lemma}[theorem]{Lemma}

\newtheorem{proposition}[theorem]{Proposition}
\newtheorem*{assumption}{Assumption}

\newtheorem*{TheoremA}{Theorem A}
\newtheorem*{TheoremB}{Theorem B}
\newtheorem*{TheoremC}{Theorem C}
\newtheorem*{TheoremD}{Theorem D}
\newtheorem*{TheoremE}{Theorem E}

\theoremstyle{definition}
\newtheorem{definition}[theorem]{Definition}
\newtheorem{example}[theorem]{Example}
\newtheorem*{Example}{Example}
\newtheorem*{notation}{Notation}
\newtheorem*{RefDiag}{Reference Diagrams}

\theoremstyle{remark}
\newtheorem{remark}[theorem]{Remark}
\newtheorem*{Remark}{Remark}

\newcommand{\leftmapsto}{\raisebox{5pt}{\rotatebox{180}{$\mapsto$}}}

\numberwithin{equation}{section}

\begin{document}

\title[Combinatorial $\mathbb{R}$-trees as generalized Cayley graphs]{Combinatorial $\mathbb{R}$-trees as\\ generalized Cayley graphs for \\fundamental groups of one-dimensional spaces}

\author{Hanspeter Fischer}

\address{Department of Mathematical Sciences\\Ball State University\\Muncie, IN 47306\\U.S.A.}

\email{fischer@math.bsu.edu}

\author{Andreas Zastrow}

\address{Institute of Mathematics, University of Gda\'nsk, ul. Wita Stwosza 57,
80-952 Gda\'nsk, Poland}

\email{zastrow@mat.ug.edu.pl}

\thanks{}

\subjclass[2000]{20F65; 20E08, 55Q52, 57M05, 55Q07}

\keywords{$\mathbb{R}$-tree, generalized Cayley graph, one-dimensional space}

\date{June 27, 2011}

\commby{}

\begin{abstract}
In their study of fundamental groups of one-dimensional path-connected compact metric
spaces, Cannon and
Conner have asked: Is there a tree-like object that might be
considered the topological Cayley graph? We answer  this question in the positive and provide a  combinatorial description of
such an object.
\end{abstract}

\mbox{\hspace{1pt}}\vspace{-.7in}

 \maketitle

\vspace{-30pt}

\tableofcontents

\vspace{-35pt}

\section{Introduction}

\noindent Fundamental groups of one-dimensional Peano continua are notoriously difficult to analyze \cite{E1, E2, E3, ADTW}.
 They are free if and only if the underlying space is locally simply-connected \cite[Theorem~2.2]{CF2}. Yet, every finitely generated subgroup of the fundamental group of a one-dimensional separable metric space is free \cite[Section~5]{CF1} and the homotopy class of every loop contains an essentially unique shortest representative (see \cite[Lemma 3.1 and Theorem 3.1]{CF2} or \cite[Theorem~3.9]{CC3}).
In light of these and related results, Cannon and
Conner have asked whether a general one-dimensional path-connected compact metric
space $X$ admits a tree-like object that might be
considered the ``topological Cayley graph'' of its fundamental
group $\pi_1(X,x)$ \cite[Question~3.9.1]{CC3}. In this article, we answer this question in the positive and provide a  combinatorial description of
such an object.

The main feature of a classical Cayley graph (for a finitely generated group) is that its vertex set bijectively corresponds to the elements of the group in such a way that the various edge-paths between two fixed vertices describe all possible representations of the difference of the corresponding group elements by words in the generators. The word length distance agrees with the natural path length metric of the Cayley graph and  the group acts by graph automorphism on the Cayley graph; it acts freely and transitively on the vertex set.

In the ``tame'' case, where the underlying space is a  one-dimensional simplicial complex, we have a free fundamental group whose Cayley graph can readily be built from the universal covering space by collapsing the lifts of a maximal subtree of the covered graph---making the vertex set of the Cayley graph the preimage of a single base vertex. The fact that the Cayley graph is a simplicial tree in this case is witness to the principle that the free group structure is fully captured by the concatenation of words and their reduction via cancellation.

\pagebreak
\begin{figure}[h!] \fbox{\parbox{4.5in}{\hspace{.2in} \includegraphics{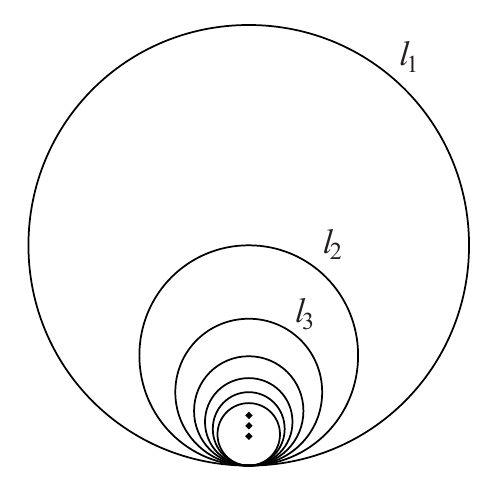} \includegraphics{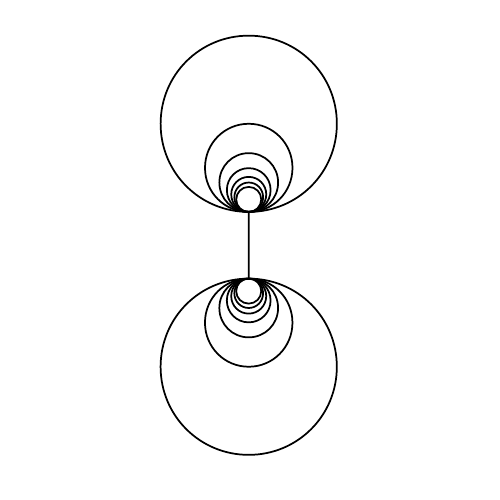} }}\caption{\label{HH} The Hawaiian Earring (left) is the one-dimensional planar set  $\mathbb{H}=\bigcup_{n\in \mathbb{N}}\{(x,y)\in \mathbb{R}^2\mid x^2+(y-1/n)^2=(1/n)^2\}$}  \end{figure}
The general situation is  more delicate.  Since  $X$ allows for the accumulation of small essential loops, we are faced with the following obstacles: (1)~The fundamental group might be uncountable; (2) there might not be a universal covering space; and (3)
collapsing a contractible subset of $X$ might drastically alter its fundamental group.
   (For example, if we collapse an arc that connects the distinguished points of two copies of the Hawaiian Earring, as depicted in Figure~\ref{HH}; see \cite[Theorem~1.2]{E2}.)

It is shown in \cite[Theorem 4.10 and Example 4.14]{FZ2} that $X$ admits a {\sl generalized} universal covering  $q:\tilde{X}\rightarrow X$ on which  $\pi_1(X,x)$ acts as the group of covering transformations, and
that $\tilde{X}$ is an $\mathbb{R}$-tree. (An $\mathbb{R}$-tree is a uniquely arcwise connected metric space in which every arc is an isometric embedding of a compact interval of the real line). We choose this $\mathbb{R}$-tree as the underlying space for our {\sl generalized} Cayley graph, keeping in mind two inevitable limitations: We must abandon the idea of using a conventional generating set, because collapsing is not an option and because we are dealing with ``nearly free''  groups that are not free on any generating set. Furthermore, there is no $\mathbb{R}$-tree metric on $\tilde{X}$ for which the action of $\pi_1(X,x)$ could possibly be by isometry. (See also Remark~\ref{noiso}.) From this point of view, the following seems to be the best possible solution to the given problem.

We give a fully combinatorial description of the $\mathbb{R}$-tree $\tilde{X}$ and its designated subset  $q^{-1}(x)=\pi_1(X,x)$ by
uniquely labeling all points with infinite sequences of finite words, which combinatorially capture the
structure of $\pi_1(X,x)$ by way of term-wise concatenation and reduction. Here, we are limited to using
 sequences of (reduced) words which specify (homotopy classes of) edge-paths through various approximating graphs for $X$, rather than the usual words whose individual letters correspond to homotopy classes of entire loops.  Arcs between two points of $q^{-1}(x)=\pi_1(X,x)$ naturally spell out word sequences that represent the difference of the corresponding group elements. We recursively assign weights to the individual letters of the words of a word sequence, in such a way that we obtain a limiting word length function which combinatorially describes the  $\mathbb{R}$-tree metric on $\tilde{X}$ as a radial metric.

In particular, we provide a combinatorial description of the fundamental group  $\pi_1(X,x)$ of a general one-dimensional path-connected compact metric space $X$ via a word calculus in which there are no relations, other than cancellation---underscoring the nearly free character of the group.

There are many situations in which $X$ is the limit of a preferred inverse system of approximating graphs, making the set-up of this paper a rather natural and systematic start of inquiry into its fundamental group. Such is the case for one-dimensional CAT(0) boundaries \cite{BH,CF}. The Sierpi\'nski carpet and the Menger universal curve, for example, arise in this way as Gromov boundaries of hyperbolic Coxeter groups \cite{B,KB}.

\section{Informal Overview of Definitions and Results}

 \begin{RefDiag} {\sl As a visual guide for the informal overview contained in this section, the reader may wish to consult the diagrams of Remarks~\ref{GroupDiagram} and \ref{TreeDiagram}.}
\end{RefDiag}

\begin{figure}[h!]\fbox{\parbox{4.5in}{\includegraphics{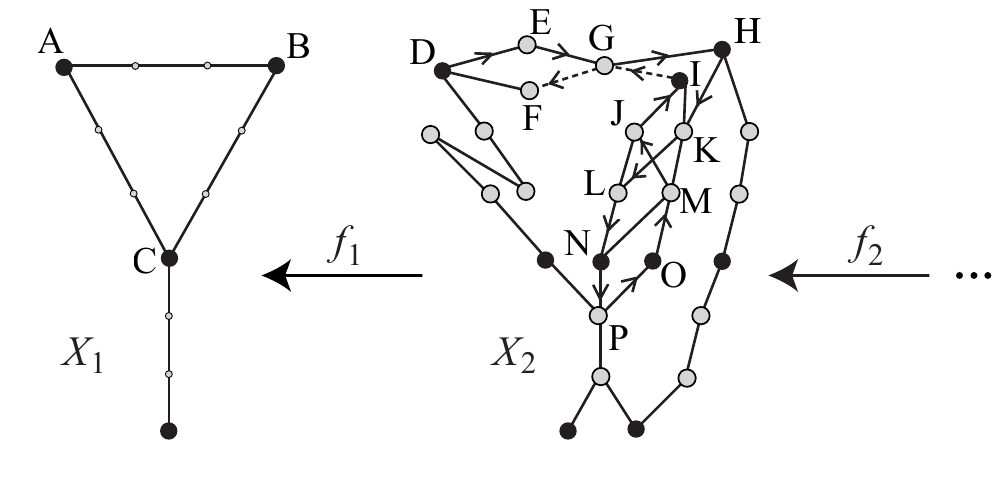}

\hspace{.25in} $\omega_1=\mbox{ABCB}$ \hspace{.335in} $\stackrel{\phi_1}{\leftmapsto}$ \hspace{.165in}  $\omega_2=\mbox{DEGHKLNPOMJI}$ \hspace{.37in} $\stackrel{\phi_2}{\leftmapsto}$ \hspace{3pt} $\cdots$

\hspace{.25in} $\omega_1=\mbox{ABCB/A}$ \hspace{.165in} $\stackrel{\phi_1}{\leftmapsto}$ \hspace{.165in}  $\omega_2=\mbox{DEGHKLNPOMJIGF}$ \hspace{.174in} $\stackrel{\phi_2}{\leftmapsto}$ \hspace{3pt} $\cdots$

 }}
\caption{\label{example}  Examples of word sequences $(\omega_1,\omega_2,\omega_3,\ldots)\in {\mathcal W}$. The smaller vertices on the left are the subdivision vertices of $X_1$. Black and gray vertices map to black or gray vertices, respectively. For example, $f_1(D)=A$, $f_1(H)=f_1(I)=B$, $f_1(N)=f_1(O)=C$. The arrows in the (undirected) graph $X_2$ indicate two possibilities for a path $\omega_2$.
 In both cases, $\omega_2$ is a reduced path, while $\omega_1$ is not.}
\end{figure}

We express the  space $X$ as the limit of an inverse sequence of finite graphs $X_n$ and bonding maps $f_n$ which map each edge of a given graph homeomorphically onto an edge of a subdivision of the previous graph. (See Figure~\ref{example} and Lemma~\ref{setup} \& Notation.)
 Edge-paths through these graphs will be recorded by words of visited vertices. Observe that the process of cancelling an adjacent inversely directed edge-pair in an edge-path generates the path homotopy classes for a given graph and that each such homotopy class contains a unique reduced representative.
 (For example, one edge-cancellation within the word $\omega_1=ABCB$ reduces it to $\omega_1'=AB$.)

 The topological bonding maps are then replaced by combinatorial word functions $\phi_n$, which naturally transform the (unreduced) words of one level into well-formed (unreduced) words of the previous level.
 All ensuing combinatorial notions are subsequently framed in terms of this combinatorial inverse limit of sets of words, denoted by $\mathcal W$, whose elements we call {\em word sequences}. (See Definition~\ref{W}.)
  The convention of suppressing adjacent repetitions of letters within a given word ensures that the words of a word sequence remain finite even if they oscillate increasingly at finer approximation stages.

Word sequences will always start at a fixed base point. Naturally, when investigating the fundamental group, word sequences will also return to the base point. (The set of returning word sequences will be denoted by $\Omega$.) When they do not, certain round-off information will have to be encoded in the ending of each word: We will signify a combinatorial end of a path ``between two vertices'' with a slash ``/'' between the last two letters of a word, as suggested in Figure~\ref{example}.
This naturally leads to a certain degree of combinatorial redundancy in the word endings of word sequences,
 similar to (but more varied than) the nonuniqueness of decimal representations (such as $0.999\ldots\doteq 1.000\ldots$), because we are approximating continuous entities by discrete objects, some of which can be approximated from different sides.
Accordingly, the symbol ``$\doteq$'' will be used to indicate that two word sequences are equal up to a combinatorially equivalent ending. We place a dot ``$\dot{\mbox{\hspace{5pt} }}$'' over an entire set of word sequences when selecting canonical representatives with respect to this equivalence relation. (Formal definitions of these concepts are given in Section~\ref{words}.)

 The elements of $\tilde{X}$ are homotopy classes of paths in $X$ that emanate from the base point $x$  and the  map $q:\tilde{X}\rightarrow X$ consists of the standard endpoint projection.
 When endowed with the correct topology, this generalized universal covering space is characterized by the usual unique lifting criterion and $\pi_1(X,x)$  acts naturally on $\tilde{X}$ as the group of covering transformations \cite{FZ2}.
 There is a natural injective map from $\tilde{X}$ into the inverse limit $\hat{X}$ of the simplicial trees which cover the finite approximating graphs of $X$.  Along with it comes a natural injective homomorphism from the fundamental group $\pi_1(X,x)$ into the first \v{C}ech homotopy group $\check{\pi}_1(X,x)$, which is the inverse limit  of the free fundamental groups of these finite graphs. (See Lemma~\ref{CurtisFort} and Remark~\ref{simply}.)

This poses the challenge of combinatorially identifying the homomorphic image of $\pi_1(X,x)$ in $\check{\pi}_1(X,x)$.
 Our solution to this problem is modeled on the work of \cite{ADTW} for the Sierpi\'nski gasket and proceeds as follows.
An element of $\check{\pi}_1(X,x)$ has a natural representation by a sequence $(g_n)_n$ of unique canonically {\sl reduced} words~$g_n$, each of which represents an entire homotopy class of edge-paths. Such a sequence  $(g_n)_n$ is not $\phi_n$-coherent (and hence not a word sequence of~${\mathcal W}$) but only $\phi_n'$-coherent, where $\phi_n'$ denotes $\phi_n$ followed by reduction. (Figure~\ref{example} shows  examples of reduced words $\omega_2$ which map to unreduced words $\omega_1$ under~$\phi_2$.)
We will denote the set of all  returning $\phi'_n$-coherent reduced sequences by $G$ and use the symbol ``$\;'\;$'' throughout when reducing words. (See Definition~\ref{prime}, Lemma~\ref{reduction} and Remark~\ref{check}.) Then for each element $g\in \pi_1(X,x)$ there is  some
 sequence $\varphi(g)=(g_n)_n\in G\cong\check{\pi}_1(X,x)$ of reduced words representing the image of $g$ in $\check{\pi}_1(X,x)$. (See Definitions~\ref{spelling} and \ref{phidef},  and Lemma~\ref{cech}.)
 If we represent an arbitrary element of $\check{\pi}_1(X,x)$ by a sequence $(g_n)_n\in G$ and project  progressively later words $g_k$ of this sequence  onto fixed lower levels $n$ without reducing them, then this process might or might not stabilize to an overall $\phi_n$-coherent word sequence $(\omega_n)_n$ of ${\mathcal W}$. If it stabilizes, at all levels, we call $(g_n)_n$ {\em locally eventually constant} and we place the symbol ``$\;\overleftarrow{\mbox{\raisebox{5pt}{\hspace{10pt}}}}\;$'' over it to denote the resulting stabilized word sequence: $\overleftarrow{(g_n)_n}=(\omega_n)_n$. (See Definition~\ref{stab}.)
Denoting by ${\mathcal G}\subseteq G\cong\check{\pi}_1(X,x)$  the set of all elements  of $G$ which stabilize in this sense, it turns out that $\varphi(\pi_1(X,x))={\mathcal G}=\Omega'$. (See Lemma~\ref{setequal} and Theorem~\ref{group}.)

The stabilized state of a reduced word sequence captures the ideal degree of combinatorial reduction, leading to a combinatorial description of $\pi_1(X,x)$ in terms of word sequences which generalizes the description given in \cite{ADTW} for the fundamental group of the Sierpi\'nski gasket:

Theorem~A of Section~\ref{results} describes the fundamental group of $X$  as the combinatorially well-formed set of word sequences $\overleftarrow{\mathcal G}$ along with the combinatorially well-defined binary operation of term-wise concatenation of words, followed by reduction and restabilization.

Similarly, every element of $\hat{X}$ can be represented by a non-returning sequence of reduced words. We will denote the set of all reduced $\phi_n'$-coherent sequences by $R$ and we will  denote the subset of sequences of $R$ that stabilize in the above sense by ${\mathcal R}\subseteq R$. We then combinatorially identify the image of $\tilde{X}$ in $\hat{X}$ in terms of word sequences by $\overleftarrow{\varphi(\tilde{X})}=\dot{\overleftarrow{\mathcal R}}$. (See Theorem~\ref{image}(b).)

In Theorem~B  we show that $\dot{\overleftarrow{\mathcal R}}$ is an $\mathbb{R}$-tree whose metric is radially induced by a word length function for word sequences. This word length function is based on a  recursive weighting scheme from \cite{MO}, applied to the letters of words of adjacent levels. (See Definitions~\ref{weights} and \ref{overlap}; see Definition~\ref{formal} for ``$DRC$''.) In order to correctly capture the topology of the $\mathbb{R}$-tree, however, the word sequences need to first undergo a combinatorial completion step which inserts limiting letters into their words. We will use the symbol ``$\;\overline{\mbox{\raisebox{5pt}{\hspace{10pt}}}}\;$'' for completion. (See Definition~\ref{CompletionDef} and Figure~\ref{L-picture}.)
Geometrically, the completed state of the word sequence can be generated by connecting the corresponding point of the $\mathbb{R}$-tree with an arc to the base point and reading off the resulting sequence of words in the finite approximating graphs. (See Corollary~\ref{coarc} and Example~\ref{L-space}.)

Theorem~C states that arcs in the $\mathbb{R}$-tree whose endpoints correspond to elements of $\pi_1(X,x)$ naturally spell out word sequences which represent the (completed state of the) difference of the group elements.

 Theorem~D combinatorially describes the action of the fundamental group on what can now be regarded as its generalized Cayley graph.
Finally, Theorem~E  presents a combinatorial criterion (cf.\@ Definition~\ref{multi}) for when the quotient under this action is homeomorphic to the original one-dimensional space.

\begin{Remark}
Any attempt to combinatorially describe the fundamental group of a space which allows for the accumulation of small essential loops requires some concept of infinite products that accounts for this effect.
The combinatorial description of the fundamental group of the Hawaiian Earring alone has been the subject of a number of papers \cite{CC1,E1,MM,Z1}, where essentially three different approaches have emerged:  (i)~studying the inverse limit of free groups which contains the given fundamental group as a subgroup; (ii)~accommodating products of infinite linear order type; or (iii) using infinite sequences of well-formed finite words (i.e., word sequences) along with well-defined combinatorial multiplication rules.
Roughly speaking, infinite products arise as limiting objects from word sequences and, in turn, word sequences can be obtained from infinite products via successively finer approximations.
While for the Hawaiian Earring the majority of authors seem to prefer the infinite product approach, all advances into combinatorial descriptions of fundamental groups of spaces with more than one accumulation point of small essential loops use, in principle, word sequences \cite{ADTW,DS,FZ1,Z2}.
\end{Remark}

\section{General Setup: Word Sequences}\label{words}

\begin{assumption} Let $X$ be a one-dimensional path-connected compact metric space.
\end{assumption}

It is well-known that $X$ can be expressed as the limit of an inverse sequence of finite graphs \cite[Theorem~1]{MS}, i.e., of finite connected one-dimensional simplicial complexes (without looping edges or multiple edges between the same two vertices).  Moreover, given any inverse sequence of finite graphs and continuous maps whose limit is $X$, there is a systematic procedure for improving the representation:

\begin{lemma}[\cite{R}]\label{setup} There is an inverse sequence $X_1\stackrel{f_1}{\longleftarrow} X_2 \stackrel{f_2}{\longleftarrow} X_3 \stackrel{f_3}{\longleftarrow} \cdots$ of \linebreak  finite connected one-dimensional simplicial complexes $X_n$ and continuous surjections $f_n:X_{n+1}\rightarrow X_n$, along with subdivisions $X^\ast_n$ of $X_n$, such that the following hold:
\begin{itemize}
\item[(a)] $\displaystyle X=\lim_{\longleftarrow} \left(X_1\stackrel{f_1}{\longleftarrow} X_2 \stackrel{f_2}{\longleftarrow} X_3 \stackrel{f_3}{\longleftarrow} \cdots\right)$.
    \item[(b)] Every edge of $X_n$ is evenly subdivided into
    the same number of edges of $X^\ast_n$. $($This number, which is assumed to be greater than 1, depends on $n$.$)$\vspace{2pt}

\item[(c)]   $f_n:X_{n+1}\rightarrow X^\ast_n$
maps every edge of $X_{n+1}$ linearly onto an edge of $X^\ast_n$.
\end{itemize}
\end{lemma}

\begin{notation}\sl We will fix a description of $X$ as given in Lemma~\ref{setup}. Throughout the paper,
elements of (and functions into) a limit of an inverse sequence will be denoted as coherent sequences of points of (and functions into) the individual terms.
\end{notation}

\begin{proof} Lemma~\ref{setup} follows from the proof of \cite[Theorem~2]{R}, upon adding further subdivision points in the inductive step \cite[Theorem~1]{R}  to ensure that (b) holds.
\end{proof}

\begin{definition}[Word sequences: $\Omega\subseteq {\mathcal W}$]\label{W}
Let $V_n$ and $E_n$ denote the vertex set and the (undirected) edge set of $X_n$, respectively. We may assume that $V_i\cap V_j=\emptyset$ for all $i\not=j$.
Let ${\mathcal P}_n$ denote the set of all non-stagnating words $v_1v_2\cdots v_k$  over the alphabet $V_n$  (i.e., $v_i\not=v_{i+1}$ for all $i=1,2,\cdots, k-1$) which describe edge-paths in $X_n$.
For convenience, we also include the empty word in ${\mathcal P}_n$.

For each word $v_1v_2\cdots v_{k}v_{k+1}\in {\mathcal P}_n$, we also form a word $v_1v_2\cdots v_{k}/v_{k+1}$ in which we symbolically separate the last letter. (We think of this new word as an edge-path which passes vertex $v_k$, but does not quite reach vertex $v_{k+1}$.) We will write
$v_1v_2\cdots v_{k}/\ast$ when discussing issues pertaining to both types of words, referring to $v_1, v_2, \dots, v_{k}$ as the {\em proper} letters.
We define \[{\mathcal P}^+_n=\{v_1v_2\cdots v_{k}/v_{k+1}\mid v_1v_2\cdots v_{k}v_{k+1}\in {\mathcal P}_n\}\] and let $\phi_n:{\mathcal P}_{n+1}\cup {\mathcal P}^+_{n+1}\rightarrow {\mathcal P}_n\cup {\mathcal P}^+_n$ denote the natural projection function, formally described in Definition~\ref{formal} below.

 Fix a base point $x=(x_n)_n\in X$ such that $x_n\in V_n$ for all $n$.
Let ${\mathcal W}_n$ be the set of all words in ${\mathcal P}_n\cup {\mathcal P}^+_n$ that start with $x_n$ and let $\Omega_n$ be the set of all words in ${\mathcal P}_n$ that start and end with $x_n$.
We define the set ${\mathcal W}$ of {\em word sequences} by
\[{\mathcal W}=\lim_{\longleftarrow} \left({\mathcal W}_1\stackrel{\phi_1}{\longleftarrow} {\mathcal W}_2
\stackrel{\phi_2}{\longleftarrow}{ \mathcal W}_3 \stackrel{\phi_3}{\longleftarrow}\cdots\right)\]
along with its subset
\[\Omega=\lim_{\longleftarrow}\left(\Omega_1\stackrel{\phi_1}{\longleftarrow} \Omega_2 \stackrel{\phi_2}{\longleftarrow} \Omega_3 \stackrel{\phi_3}{\longleftarrow} \cdots\right).\]
\end{definition}

\begin{remark}
In Section~\ref{results}, we will represent the identity element of $\pi_1(X,x)$ by the word sequence $(\omega_n)_n\in {\mathcal W}$ with $\omega_n=x_n$ for all $n$. Accordingly, the word length function of Definition~\ref{weights} will assign a value of zero to this word sequence.
\end{remark}

\begin{definition}[Delete-Replace-Compress: ``$DRC$'' and $\phi_n$]  \label{formal}
For a given word  $\omega_{n+1}=v_1v_2\cdots v_k\in {\mathcal P}_{n+1}$, we let $DRC_n(\omega_{n+1})\in {\mathcal P}_n$ be the word obtained from $\omega_{n+1}$ by first {\sl deleting} every letter $v$ from $\omega_{n+1}$ for which $f_n(v)\not\in V_n$, next {\sl replacing} every remaining letter $v$ by $f_n(v)$, and finally {\sl compressing} any resulting maximal stagnating subwords of the form $uu\cdots u$  into one letter $u$.

We then define $\phi_n:{\mathcal P}_{n+1}\cup {\mathcal P}^+_{n+1}\rightarrow {\mathcal P}_n\cup {\mathcal P}^+_n$  as follows:

\begin{itemize}
\item[(1)] Suppose  $\omega_{n+1}=v_1v_2\cdots v_k\in {\mathcal P}_{n+1}$.
If $f_n(v_k)\in V_n$ or if $DRC_n(\omega_{n+1})$ is the empty word, then we define $\phi_n(\omega_{n+1})=DRC_n(\omega_{n+1})\in {\mathcal P}_n$; otherwise we consider $j=\max\{i\mid 1\leqslant i \leqslant k-1, f_n(v_i)\in V_n\}$ and define $\phi_n(\omega_{n+1})=DRC_n(\omega_{n+1})/u\in {\mathcal P}^+_n$, where $\{f_n(v_j),u\}\in E_n$ is the edge containing $f_n(v_{j+1})$.

\item[(2)]
     Suppose $\omega_{n+1}=v_1\cdots v_{k}/v_{k+1}\in {\mathcal P}^+_{n+1}$. If $DRC_n(v_1v_2\cdots v_k)$ is the empty word, then we define $\phi_n(\omega_{n+1})$ to be the empty word; otherwise we consider $j=\max\{i\mid 1\leqslant i \leqslant k, f_n(v_i)\in V_n\}$ and define $\phi_n(\omega_{n+1})=DRC_n(v_1v_2\cdots v_{k})/u$,  where $\{f_n(v_j),u\}\in E_n$ contains $f_n(v_{j+1})$.
 \end{itemize}
 \end{definition}

 \begin{remark} We always have  $\phi_n(v_1v_2\cdots v_k/\ast)=DRC_n(v_1v_2\cdots v_k)/\ast$.

 \end{remark}
\begin{remark} By definition,
 $\phi_n|_{\Omega_{n+1}}=DRC_n$.
\end{remark}
\begin{definition}[Terminating type] We categorize word sequences $(\omega_n)_n\in {\mathcal W}$ into two types.
\begin{itemize}
\item[(1)] {\em Terminating type:}
 there is an $N\in \mathbb{N}$ such that $\omega_n\in {\mathcal P}_n^+$ for all $n< N$ and $\omega_n\in {\mathcal P}_n$ for all  $n\geqslant N$;

\item[(2)] {\em Non-terminating type:} $\omega_n\in {\mathcal P}_n^+$ for all $n$.
\end{itemize}
\end{definition}

 \begin{remark}\label{last}
 For a word sequence $(\omega_n)_n\in {\mathcal W}$ of terminating type, $\phi_n$ maps the last letter of $\omega_{n+1}$ to the last letter of $\omega_n$ for all $n\geqslant N$.
\end{remark}

\begin{remark} Every  $(\omega_n)_n\in \Omega$ is of terminating type (with $N=1$).
\end{remark}

We now define a word sequence analog to ``$0.999\ldots\doteq 1.000\ldots$''.

\begin{definition}[Equivalence: $(\xi_n)_n \doteq(\omega_n)_n$. Terminating representatives: $\dot {\mathcal S}\;$]\label{comp}
Let $(\omega_n)_n\in {\mathcal W}$ be a word sequence of terminating type.
We call a word sequence $(\xi_n)_n\in {\mathcal W}$ of non-terminating type {\em formally equivalent} to $(\omega_n)_n$ and we write $(\xi_n)_n \doteq(\omega_n)_n$, if there is an index $N$ such that $\omega_n=v_{n,1}v_{n,2}\cdots  v_{n,{m_n}}$ for all $n\geqslant N$, and
either $\xi_n=v_{n,1}v_{n,2}\cdots v_{n,{m_n-1}}/v_{n,{m_n}}$  for all $n\geqslant N$,
 or   $\xi_n=v_{n,1}v_{n,2}\cdots v_{n,{m_n}}/v_{n,{m_n}+1}$  for all $n\geqslant N$ and some $v_{n,{m_n}+1}$.
 We denote the induced equivalence relation on $\mathcal W$ also by the symbol $\doteq$.
 Given  $\mathcal S\subseteq \mathcal W$, we denote by $\dot {\mathcal S}$ the set of word sequences obtained via replacing every element of $\mathcal S$ by a formally equivalent element from $\mathcal W$ of terminating type, whenever possible.
\end{definition}

\begin{remark} Note that in Definition~\ref{comp}, we might not be able to choose $N$ so that $\xi_n=\omega_n$ for all $n<N$. Indeed, the relationship between the words $\xi_n$ and $\omega_n$ might be reversed for some $n<N$ when compared to $n\geqslant N$. Specifically, we may have $\xi_n=v_{n,1}v_{n,2}\cdots v_{n,{m_n-1}}/v_{n,{m_n}}$  for all~$n$ while $\omega_n=v_{n,1}v_{n,2}\cdots  v_{n,{m_n}}$ for all $n\geqslant N$ and
  $\omega_k=v_{k,1}v_{k,2}\cdots  v_{k,{m_k}-1}$ for one or more $k<N$. This will happen when for some $1\leqslant i< m_N-1$, each of the three words $v_{N,i}$, $v_{N,m_N}$, and $v_{N,i}v_{N,i+1}\cdots v_{N,m_N}$ gets mapped to the letter $v_{k,m_k-1}$ by $\phi_k\circ\phi_{k+1}\circ\cdots\circ\phi_{N-1}$.
\end{remark}

\begin{remark}\label{classes}
If a formal equivalence class of $\mathcal W$ contains more than one element, then it contains exactly one word sequence of terminating type and possibly uncountably many word sequences of non-terminating type.
\end{remark}

\section{Combinatorial Notions and Definitions} \label{combinatorics}

{\sl The definitions of this section are solely in terms of the functions $\phi_n$.}

\begin{definition}[Concatenation: $\omega_n\xi_n$]
For two words $\omega_n=v_1v_2\cdots v_k\in \Omega_n$ and\linebreak $\xi_n=u_1u_2\cdots u_s/\ast \in {\mathcal W}_n$ we define $\omega_n\xi_n=v_1v_2\cdots v_{k-1}u_1u_2\cdots u_s/\ast\in {\mathcal W}_n$.
\end{definition}

\begin{definition}[Reduction: $(\omega_n)'_n$, ${\mathcal W}'\subseteq R$, $\Omega'\subseteq G$]\label{prime}
The {\em reduction} $\omega'_n\in {\mathcal W}_n$ of a given word  $\omega_n\in {\mathcal W}_n$ is obtained by repeatedly replacing substrings of $\omega_n$ of the form ``$uvu$'' and ``$uv/u$'' by ``$u$'' and ``$u/v$'', respectively, until this is no longer possible.  We will call $\omega_n$ {\em  reduced} if $\omega'_n=\omega_n$. Consider the set ${\mathcal W}_n'=\{\omega_n'\mid \omega_n\in {\mathcal W}_n\}=\{\omega_n\in {\mathcal W}_n \mid \omega_n \mbox{ is reduced}\}$ of all reduced words in ${\mathcal W}_n$
 and let $\phi'_n:{\mathcal W}_{n+1}'\rightarrow {\mathcal W}_n'$ be the function given by $\phi'_n(\omega'_{n+1})=\phi_n(\omega'_{n+1})'$. We define the set $R$ by  \[R=\lim_{\longleftarrow} \left({\mathcal W}_1'\stackrel{\phi'_1}{\longleftarrow}  {\mathcal W}_2'\stackrel{\phi'_2}{\longleftarrow} {\mathcal W}_3'\stackrel{\phi'_3}{\longleftarrow}\cdots \right).\] We also define a subset $G\subseteq R$ by considering the set $\Omega_n'$ of all reduced words in~$\Omega_n$, i.e., $\Omega_n'=\{\omega_n'\mid \omega_n\in \Omega_n\}\subseteq {\mathcal P}_n$, and setting \[ G=\lim_{\longleftarrow} \left(\Omega_1'\stackrel{\phi'_1}{\longleftarrow} \Omega_2'\stackrel{\phi'_2}{\longleftarrow} \Omega_3'\stackrel{\phi'_3}{\longleftarrow}\cdots \right).\]
 Moreover, for a word sequence $(\omega_n)_n\in{\mathcal W}$, we define $(\omega_n)'_n=(\omega'_n)_n$, and  for a subset ${\mathcal S}\subseteq {\mathcal W}$, we define ${\mathcal S}'=\{(\omega_n)_n'\mid (\omega_n)_n\in {\mathcal S}\}$. Since $\phi_n(\omega_{n+1})'=\phi_n(\omega_{n+1}')'$ for all $(\omega_n)_n\in {\mathcal W}$ and all $n$, we have  ${\mathcal W}'\subseteq R $ and $\Omega'\subseteq G$.
\end{definition}

\begin{remark} By Lemma~\ref{reduction} below, reduction is well-defined.
\end{remark}

\begin{remark} In general, $R \nsubseteq {\mathcal W}$ and $G \nsubseteq \Omega$, because the sequences of $R$ and $\Omega$ are $\phi'_n$-coherent rather than $\phi_n$-coherent.
Moreover,  ${\mathcal W}'\subsetneqq R$ and $\Omega'\subsetneqq G$, in general. This is best illustrated by considering the sequence of reduced words that describe the commutators $l_1l_2l_1^{-1}l_2^{-1}l_1l_3l_1^{-1}l_3^{-1}\cdots l_1l_nl_1^{-1}l_n^{-1}$ in the approximating graphs of an appropriately chosen inverse sequence whose limit is the Hawaiian Earring depicted in Figure~\ref{HH}. This sequence lies in $G$ but neither in ${\mathcal W}$ nor in ${\mathcal W}'$.
\end{remark}

\begin{remark}\label{check}
Each $\Omega_n'$ forms a free group under the operation $\omega_n\ast \xi_n=(\omega_n\xi_n)'$.
Every $\phi'_n:\Omega_{n+1}'\rightarrow \Omega_n'$ is a homomorphism and the group $G$ is naturally isomorphic to the first \v{C}ech homotopy group $\check{\pi}_1(X,x)$. (See Lemma~\ref{cech} below.)
\end{remark}

\begin{definition}[Stabilization: $\overleftarrow{(r_n)_n}$, ${\mathcal R}\subseteq R$, ${\mathcal G}\subseteq G$]\label{stab} We will call a sequence  $(r_n)_n\in R$  {\em locally eventually constant} if for every fixed level $n\geqslant 1$ the sequence $(\phi_n\circ \phi_{n+1}\circ \cdots \circ \phi_{k-1}(r_k))_{k>n}$  of (unreduced) words in ${\mathcal W}_n$ is eventually constant.\footnote{We adapt the terminology ``locally eventually constant'' from \cite{MM}.} We put \[{\mathcal R}=\{(r_n)_n\in R\mid (r_n)_n \mbox{ is locally eventually constant}\},\]
\[{\mathcal G}=\{(g_n)_n\in G\mid (g_n)_n \mbox{ is locally eventually constant}\}.\] For $(r_n)_n\in  \mathcal R$, let $\omega_n=\phi_n\circ \phi_{n+1}\circ \cdots \circ \phi_{k-1}(r_k)$, for sufficiently large $k$, and define $\overleftarrow{(r_n)_n}=(\omega_n)_n\in {\mathcal W}$. We call $\overleftarrow{(r_n)_n}$ the  {\em stabilization} of $(r_n)_n$.  Finally, we define \[\overleftarrow{\mathcal R}=\{\overleftarrow{(r_n)_n}\mid (r_n)_n \in {\mathcal R}\}\subseteq {\mathcal W},\]
\[\overleftarrow{\mathcal G}=\{\overleftarrow{(g_n)_n}\mid (g_n)_n\in \mathcal G\}\subseteq \Omega.\]
\end{definition}

\begin{remark}\label{representation} The (reduced) locally eventually constant sequences $\mathcal R\subseteq R$ naturally correspond to
  the (unreduced) stabilized word sequences  $\overleftarrow{\mathcal R}\subseteq {\mathcal W}$ because \[\left(\overleftarrow{(r_n)_n}\right)'=(r_n)_n,\]  which follows from the fact that $\phi_n\circ \phi_{n+1}\circ \cdots \circ \phi_{k-1}(r_k)'=\phi'_n\circ \phi'_{n+1}\circ \cdots \circ \phi'_{k-1}(r_k)$.
   That is, we have the following bijection:
   \[
\xymatrix{R & \; {\mathcal R} \ar@{_{(}->}[l] \ar@{.>}@/^1pc/[rr]^{``\overleftarrow{\;\;\;}"}  & bij. & \ar@{.>}@/^1pc/[ll]^{'} \overleftarrow{{\mathcal R}}  \ar@{^{(}->}[r] & {\mathcal W}}
\]

\end{remark}

\begin{remark}\label{Rsetequal} By Lemma~\ref{setequal} below, we have ${\mathcal R}={\mathcal W}'$ and ${\mathcal G}=\Omega'$.
\end{remark}

Completion inserts limiting letters into the words of a word sequence:

\begin{definition}[Completion: $\overline{(\omega_n)_n}\;$] \label{CompletionDef} Given a word sequence $(\omega_n)_n\in {\mathcal W}$, we define its {\em completion} $\overline{(\omega_n)_n}\in {\mathcal W}$ based on the following modification of  $DRC$:

For $k>n+1$ and any word $v_1v_2\cdots v_m\in {\mathcal P}_k$, we let
$drc^k_n(v_1v_2\cdots v_m)\in {\mathcal P}_n$ be the word obtained from $v_1v_2\cdots v_m\in {\mathcal P}_k$ in three steps: first delete every letter $v$ from $v_1v_2\cdots v_m$ for which  $\phi_n\circ \phi_{n+1}\circ\cdots\circ \phi_{k-1}(v)$ is the empty word, unless there is a (unique\footnote{Here we need $k>n+1$, rather than $k>n$, because $X^\ast_n$ might only halve the edges of~$X_n$.}) letter $u$ with $uv\in {\mathcal P}_k$ such that $\phi_n\circ \phi_{n+1}\circ\cdots\circ \phi_{k-1}(u)$ is not the empty word; then replace every remaining
letter $v$ by the letter
$\phi_n\circ \phi_{n+1}\circ\cdots\circ \phi_{k-1}(v)$ or the  letter $\phi_n\circ \phi_{n+1}\circ\cdots\circ \phi_{k-1}(u)$, respectively; finally compress the resulting maximal stagnating subwords into one letter as before.

For each $n$, express $\omega_n=v_{n,1}v_{n,2}\cdots v_{n,{m_n}}/\ast$. Now fix $n$. As $k$ increases, the words  $drc^k_n(v_{k,1}\cdots v_{k,{m_k}})$ are eventually constant (see Lemma~\ref{completion}); say for $k\geqslant K$. For $k\geqslant K>n+1$, let $j_k$ be the maximal index for which $drc_n^k$ does not delete the letter $v_{k,j_k}$ from the word $v_{k,1}\cdots v_{k,{m_k}}$. If  $j_k<m_k$ for some $k\geqslant K$, then
 $\omega_n=v_{n,1}v_{n,2}\cdots v_{n,{m_n}}/v_{n,{m_n+1}}$ for some $v_{n,{m_n+1}}$ and the word   $drc^k_n(v_{k,1}v_{k,2}\cdots v_{k,{m_k}})$  ends either in the letter  $v_{n,{m_n+1}}$ or in the letter $v_{n,{m_n}}$ and, accordingly, we
   put $\tau_n=drc^k_n(v_{k,1}v_{k,2}\cdots v_{k,{m_k}})/v_{n,{m_n}}$ or $\tau_n=drc^k_n(v_{k,1}v_{k,2}\cdots v_{k,{m_k}})/v_{n,{m_n+1}}$. If  $j_k=m_k$ for all $k\geqslant K$, then we  define $\tau_n=drc^k_n(v_{k,1}v_{k,2}\cdots v_{k,{m_k}})$.
    Finally, we define $\overline{(\omega_n)_n}=(\tau_n)_n$.
\end{definition}

\begin{remark} By Lemma~\ref{completion}, the completion of a word sequence
  is well-defined.  Moreover, if $(\omega_n)_n\in \Omega$ then $\overline{(\omega_n)_n}\in\Omega$.
\end{remark}

\begin{remark}\label{insertions}
   While the process of completion inserts limiting letters into the words of a word sequence, it might also drop one letter at the  end of some of the words. Based on  the definition of $drc_n^k$, some of the proper letters ``$u$'' in the words of a word sequence might get replaced by strings of the form ``$uv_1uv_2u\cdots v_s u$'', while the ending of a word can change in one of the following three ways:  $(\dots u/v)\mapsto (\dots uv)$, $(\dots u/v)\mapsto (\dots uv/u)$, $(\dots u/v)\mapsto (\dots u)$. In particular, if $(\omega_n)_n$ is of terminating type then so is
$\overline{(\omega_n)_n}$, but not necessarily vice versa.
\end{remark}

\begin{remark} For $(\omega_n)_n\in {\mathcal W}$, in general, $\overline{(\omega_n)_n}\not\in \overleftarrow{\mathcal R}$ (cf.\@ Example~\ref{L-space}).
\end{remark}

\begin{remark} In Lemma~\ref{period}(a), we will show the following correspondence, which
improves upon Remark~\ref{representation} for returning word sequences:

\hspace{-18pt} \parbox{8in}{
\[
\xymatrix{ \overline{\overleftarrow{\mathcal G}} \ar@{.>}[dr]_{'}^{ bijections} & &  \overleftarrow{\mathcal G}\; \ar@{.>}_{``\overline{\;\;\;}"}[ll] \\
& \; {\mathcal G} \ar@{.>}_{``\overleftarrow{\;\;\;}"}[ur]}
\]}

\end{remark}

\begin{definition}[Dynamic word length: $\|(\omega_n)_n\|\;$] \label{weights} For a fixed word sequence $(\omega_n)_n\in {\mathcal W}$,
 we recursively assign weights to the letters of the words $\omega_n$ as follows.

To the letters $v_1,v_2,\ldots, v_s$ of the first word $\omega_1=v_1v_2\cdots v_s/\ast$ we assign the weights $\frac{1}{2},\frac{1}{2^2},\cdots,\frac{1}{2^s}$, respectively. (For words of the form $v_1v_2\cdots v_s/v_{s+1}$, we never assign any weight to the letter $v_{s+1}$.) Assuming that the letters $v_1,v_2,\ldots, v_k$ of the word $\omega_n=v_1v_2\cdots v_k/\ast$ have been assigned the weights $a_1,a_2,\cdots, a_k$, respectively, we assign weights $b_1, b_2, \cdots , b_m$ to the letters $u_1,u_2,\ldots, u_m$ of the word $\omega_{n+1}=u_1u_2\cdots u_m/\ast$ as follows. Since $DRC_n(u_1u_2\cdots u_m)=v_1v_2\cdots v_k$,  we may cut the word $\omega_{n+1}$ into substrings in such a way that $i_1$ is the maximal index with $DRC_n(u_1u_2\cdots u_{i_1})=v_1$ and, inductively, $i_{t+1}$ is the maximal index with $DRC_n(u_{i_t+1}u_{i_t+2}\cdots u_{i_{t+1}})=v_{t+1}$, the last index being $i_k=m$:
\[\xymatrix@1@=0pt@R=20pt{v_1& &v_2 & &v_3 & &\cdots \mbox{\hspace{1pt} } & v_k\\
u_1\ar@{|-{>}}[u]_{DRC_n}&\cdots u_{i_1}\Big|  & u_{i_1+1}\ar@{|-{>}}[u]_{DRC_n}&  \cdots u_{i_2} \Big|  & u_{i_2+1}\ar@{|-{>}}[u]_{DRC_n} & \cdots u_{i_3} \Big| &  \cdots \Big| & u_{i_{k-1}+1} \ar@{|-{>}}[u]_{DRC_n}&\cdots u_m
}
\]
We then define the weights $b_1, b_2, \cdots, b_m$ by \[\frac{a_1}{2},\frac{a_1}{2^2},\frac{a_1}{2^3}, \cdots ,\frac{a_1}{2^{i_1}}\Big| \frac{a_1}{2^{i_1}}+\frac{a_2}{2},\frac{a_2}{2^2},\frac{a_2}{2^3},\cdots, \frac{a_2}{2^{i_2-i_1}}\Big|  \frac{a_2}{2^{i_2-i_1}}+\frac{a_3}{2},\frac{a_3}{2^2},\frac{a_3}{2^3},\cdots,\frac{a_3}{2^{i_3-i_2}}\Big| \cdots \] \[\cdots \Big| \frac{a_{k-1}}{2^{i_{k-1}-i_{k-2}}}+\frac{a_k}{2},\frac{a_k}{2^2},\cdots, \frac{a_k}{2^{m-i_{k-1}}}.\]
 (Notice the carryover after each subdivision.) While $a_i$ is the weight of the $i^\text{th}$ letter of the $n^\text{th}$ word $\omega_n=v_1v_2\cdots v_k/\ast$ of the word sequence $(\omega_n)_n$, we will abuse notation and simply denote $a_i$ by $|v_i|$ whenever it is clear from context what we mean. We define the  length of the word $\omega_n=v_1v_2\cdots v_k/\ast$ as the sum of the weights of its proper letters: $|\omega_n|=|v_1|+|v_2|+\cdots +|v_k|$.  The lengths $|\omega_n|$ decrease with increasing $n$ so that we may define  the {\em length} of the entire word sequence $(\omega_n)_n$ by \[\displaystyle \|(\omega_n)_n\|=\lim_{n\rightarrow \infty} |\omega_n|.\]
\end{definition}

Next, we define to the concept of stable initial match as the maximal sub-word sequence of two word sequences:

\begin{definition}[Stable initial match: $(\omega_n)_n\Cap (\xi_n)_n\;$]\label{overlap}
For two word sequences $(\omega_n)_n,(\xi_n)_n\in {\mathcal W}$, we denote by $\omega_k\cap \xi_k$ the maximal matching consecutive initial substring of letters  of the two words $\omega_k, \xi_k\in {\mathcal W}_k$, including any letters that might come after the symbol ``/'', where we separate the last two letters of $\omega_k\cap \xi_k$ by the symbol ``/'' if they are so separated in the shorter of the two words $\omega_k$ and $\xi_k$. For $n<k_1<k_2$, the word $\phi_n\circ\phi_{n+1}\circ \cdots \circ \phi_{k_2-1}(\omega_{k_2}\cap \xi_{k_2})$ is an initial substring of
$\phi_n\circ\phi_{n+1}\circ \cdots \circ\phi_{k_1-1}(\omega_{k_1}\cap \xi_{k_1})$. Hence, with increasing~$k$, $\phi_n\circ\phi_{n+1}\circ \cdots\circ \phi_{k-1}(\omega_{k}\cap \xi_{k})$ is eventually constant; say it eventually equals $\tau_n$. We define the {\em stable initial match} of $(\omega_n)_n$ and $(\xi_n)_n$  by $(\omega_n)_n\Cap (\xi_n)_n=(\tau_n)_n\in {\mathcal W}$.
\end{definition}

  \begin{Example} We have  $(v_1v_2/v_3)\cap (v_1v_2v_3)=v_1v_2/v_3$, $(v_1v_2/v_3)\cap (v_1v_2/v_4)=(v_1v_2)$, $(v_1v_2)\cap( v_1v_2/v_3)=(v_1v_2)$, $(v_1v_2/v_3)\cap (v_1v_2v_3/v_4)=(v_1v_2/v_3)$. \qed
\end{Example}

While every letter of a given level potentially splits into multiple preimage letters at subsequent levels, its multiplicity may be essentially bounded:

\begin{definition}[Essential multiplicity]\label{multi} Fix $v\in V_n$. For each $k>n$ consider the set $V_k(v)=\{u\in V_k\mid \phi_n\circ \phi_{n+1}\circ\cdots\circ \phi_{k-1}(u)=v\}$.
  For $u_1, u_2\in V_k(v)$, we write $u_1\stackrel{v}{\sim} u_2$ if there is a word $\omega_k\in {\mathcal P}_k$ whose first letter is $u_1$ and whose last letter is $u_2$, such that the word $\phi_n\circ\phi_{n+1}\circ\cdots\circ\phi_{k-1}(\omega_k)$ consists of the single letter $v$. Let $c_k(v)$ denote the number of  $\stackrel{v}{\sim}$-equivalence classes in $V_k(v)$. The numbers $c_k(v)$ increase with $k$ and we call $\displaystyle \lim_{k\rightarrow \infty} c_k(v)$  the {\em  essential multiplicity} of $v$.
\end{definition}

\begin{Example}
In Figure~\ref{example} above, we have $|V_2(C)|=4$, $N\stackrel{C}{\sim} O$ and $c_2(C)=2$. \qed
\end{Example}

\section{Statements of Results (Theorems A--E)}\label{results}

\begin{TheoremA} The word sequences of $\overleftarrow{\mathcal G}$  form a group under the binary operation given by $(\omega_n)_n \ast (\xi_n)_n=\overleftarrow{(\omega_n\xi_n)'}$, and the group $\overleftarrow{\mathcal  G}$ is isomorphic to $\pi_1(X,x)$.
\end{TheoremA}

\noindent {\em Proof.} This theorem will be proved as Theorem~\ref{group} below.

\begin{TheoremB} For word sequences $(\omega_n)_n, (\xi_n)_n\in \overleftarrow{\mathcal R}$, define \[\rho((\omega_n)_n,(\xi_n)_n)=\Big\|\overline{(\omega_n)_n}\Big\|+\Big\|\overline{(\xi_n)_n}\Big\|-2\Big\|\overline{(\omega_n)_n}\Cap \overline{(\xi_n)_n}\Big\|.\] Then $\rho$ is a pseudo metric on $\overleftarrow{\mathcal R}$ with $\rho((\omega_n)_n,(\xi_n)_n)=0 \Leftrightarrow (\omega_n)_n\doteq (\xi_n)_n$.  Moreover, the resulting metric space $(\dot{\overleftarrow{\mathcal R}},\rho)$ is an $\mathbb{R}$-tree.
\end{TheoremB}

\noindent {\em Proof.} This theorem will be proved as Corollary~\ref{pseudo} below.

\begin{TheoremC}
For $(\omega_n)_n,(\xi_n)_n\in \overleftarrow{\mathcal G}\cong \pi_1(X,x)$, the arc of the $\mathbb{R}$-tree $\dot{\overleftarrow{\mathcal R}}$ from $(\omega_n)_n$ to $(\xi_n)_n$ naturally spells out the word sequence  $\overline{(\omega_n)_n^{-1}\ast (\xi_n)_n}$.
\end{TheoremC}

\noindent {\em Proof.}
This theorem will be proved as Corollary~\ref{differences} below.

\begin{TheoremD}
The group $\overleftarrow{\mathcal G}\cong \pi_1(X,x)$ acts freely and by homeomorphism on the $\mathbb{R}$-tree $\dot{\overleftarrow{\mathcal R}}$ via its natural action $(\omega_n)_n.(\xi_n)_n=\overleftarrow{(\omega_n \xi_n)'_n}$.
\end{TheoremD}

\noindent {\em Proof.} This theorem will be proved as Corollary~\ref{ThmD} below.

\begin{TheoremE}
 If the essential multiplicity of every letter is finite, which happens precisely when $X$ is locally path-connected, then $\dot{\overleftarrow{\mathcal R}}/\overleftarrow{\mathcal G}$ is homeomorphic to $X$.
\end{TheoremE}

\noindent {\em Proof.} This theorem will be proved as Theorem~\ref{lc} below.

\begin{remark}\label{noiso}
In general, the action of $\overleftarrow{\mathcal G}$ on $\dot{\overleftarrow{\mathcal R}}$ is not by isometry. In fact, when $X$ is the Hawaiian Earring, then there is no $\mathbb{R}$-tree metric for the topology of $\dot{\overleftarrow{\mathcal R}}$ that would render the action of $\overleftarrow{\mathcal G}$ as isometries \cite[Example~4.15]{FZ2}.
\end{remark}

\section{Proofs}

\begin{lemma}\label{reduction}
The reduction $\omega'_n$ of a given word $\omega_n\in {\mathcal W}_n$ is well-defined.
\end{lemma}

\begin{proof} If $\omega_n\in {\mathcal P}_n$, then
  $\omega'_n$ corresponds to the unique shortest representative for the homotopy class of edge-paths in $X_n$ which contains the edge-path tracing out the word $\omega_n$. The same argument can be made for $\omega_n=v_1v_2\cdots v_k/u\in {\mathcal P}^+_n$, if we temporarily allow ourselves to once subdivide the edge $\{v_k,u\}\in E_n$.
\end{proof}

\begin{lemma}\label{setequal}
We have  $\Omega'={\mathcal G}$ and ${\mathcal W}'={\mathcal R}$.
\end{lemma}

\begin{proof}
First, let $(\omega_n)_n\in \Omega$ be given. We wish to show that $(\omega_n)_n'\in {\mathcal G}$.
Observe that for every $n\leqslant k$,  the word $\phi_n\circ \phi_{n+1}\circ\cdots\circ \phi_{k-1}(\omega_k')$, when regarded as a finite sequence, is a subsequence of $\phi_n\circ \phi_{n+1}\circ\cdots\circ \phi_{k-1}(\omega_k)=\omega_n$. Moreover, $\phi_n\circ \phi_{n+1}\circ\cdots\circ \phi_{k-1}(\omega_k')$ is a subsequence of $\phi_n\circ \phi_{n+1}\circ\cdots\circ \phi_{k-1}\circ\phi_k(\omega_{k+1}')$, which is in turn a subsequence of
$\phi_n\circ \phi_{n+1}\circ\cdots\circ \phi_{k-1}\circ\phi_k\circ \phi_{k+1}(\omega_{k+2}')$, etc., all of which are  subsequences of $\omega_n$ by the above observation. Hence $(\omega_n)_n'$ is locally eventually constant and we have $(\omega_n)_n'\in {\mathcal G}$.

Next, let $(g_n)_n\in {\mathcal G}$ be given. Put $(\omega_n)_n=\overleftarrow{(g_n)_n}$. Then for every $n$ and sufficiently large $k$, $\omega_n'=\phi_n\circ \phi_{n+1}\circ \cdots\circ \phi_{k-1}(g_k)'=\phi_n'\circ \phi_{n+1}\circ \cdots\circ \phi_{k-1}(g_k)=\phi_n'\circ \phi_{n+1}'\circ \cdots\circ \phi_{k-1}'(g_k)=g_n$. Hence, $(\omega_n)_n'=(g_n)_n$ so that $(g_n)_n\in \Omega'$.

 The argument for ${\mathcal W}'={\mathcal R}$ is exactly the same, once we generalize the notion of subsequence to elements of ${\mathcal P}_n^+$ in the obvious way: $u_1u_2\cdots u_k/u_{k+1}$ is a {\em subsequence} of $v_1v_2\cdots v_m/v_{m+1}$ if $u_1u_2\cdots u_k$ is a subsequence of $v_1v_2\cdots v_m$ and $\{u_k,u_{k+1}\}=\{v_m,v_{m+1}\}$.
\end{proof}

\begin{lemma}\label{completion}
For $(\omega_n)_n\in {\mathcal W}$, the completion $\overline{(\omega_n)_n}$  is well-defined and $\overline{(\omega_n)_n}\in\mathcal W$.
\end{lemma}

\begin{proof} For every $k>n+1$, the first $m_n$ letters of the word $\omega_n=v_{n,1}v_{n,2}\cdots v_{n,{m_n}}/\ast$ record those vertices of $X_n$ that are traversed by the image under the function $f_n\circ f_{n+1}\circ\cdots\circ f_{k-1}$ of the edge-path in $X_k$, which is represented by the first $m_k$ letters of the word $\omega_k=v_{k,1}v_{k,2}\cdots v_{k,{m_k}}/\ast$ (while ignoring repeats). The word
$drc_n^k(v_{k,1}v_{k,2}\cdots v_{k,{m_k}})$ records, in addition, all vertices of $X_n$ that were narrowly missed by this image. (In the process, $drc_n^k$ may also restore some of the letters that fell victim to compression due to repetition when $v_{n,1}v_{n,2}\cdots v_{n,{m_n}}$ was formed by $DRC$ from $\omega_k$.) The larger the index $k$, the nearer the miss of the vertex. Therefore, all potential ``inserts'' in the word $\omega_n$, which  $drc_n^k$ might make for large $k$, are already determined by the word $\omega_{n+2}$. By the same token, the potential ``inserts'' in $\omega_n$ are also determined by the word $\omega_{n+3}$. However, the potential inserts determined by $\omega_{n+3}$ are a subset of the potential inserts determined by $\omega_{n+2}$. Continuing with this logic, we see that $drc_n^k(v_{k,1}v_{k,2}\cdots v_{k,{m_k}})$ is eventually constant, for sufficiently large $k$. Moreover, $DRC_n\circ drc_{n+1}^k=drc_n^k$
so  that $\overline{(\omega_n)_n}\in\mathcal W$.
\end{proof}

The following technical lemma will be needed in the buildup of the diagrams of Remarks~\ref{GroupDiagram} and \ref{TreeDiagram}. It states that completing a word sequence before reducing and restabilizing it, results in a formally equivalent word sequence and that formally equivalent word sequences have identical completions.

\begin{lemma}\label{period} Let $(\omega_n)_n, (\xi_n)_n\in {\mathcal W}$.
\begin{itemize}

\item[(a)] If  $(\omega_n)_n$ is of terminating type or if $\overline{(\omega_n)_n}$ is of non-terminating type, then $\left(\overline{(\omega_n)_n}\right)'=(\omega_n)'_n$.

\item[(b)] We always have $\overleftarrow{\left(\overline{(\omega_n)_n}\right)'}\doteq\overleftarrow{(\omega_n)'_n}$.

\item[(c)] If $(\omega_n)_n \doteq(\xi_n)_n$, then $\overline{(\omega_n)_n}=\overline{(\xi_n)_n}$.

\end{itemize}
\end{lemma}

\begin{proof} (a) If $(\omega_n)_n$ is of terminating type, then any letters that the completion process might insert into the words $\omega_n$   disappear upon reduction. (See Remark~\ref{insertions}.) If $\overline{(\omega_n)_n}$ is of non-terminating type then so is $(\omega_n)_n$ and their words have the same ending pairs, with $u/v$  switched to $v/u$ by the completion process exactly when reduction reverses this switch.

(b) By Part (a), we may assume that $(\omega_n)_n$ is of non-terminating type and that $\overline{(\omega_n)_n}$ is of terminating type.
Then $\overleftarrow{(\omega_n)'_n}$ is of non-terminating type and $\overleftarrow{\left(\overline{(\omega_n)_n}\right)'}$ is of terminating type, with all of their words identical except for the endings, which for the former is always of the form $u/v$ where the latter will eventually feature the (matching) letters $uv$ or $u$, instead. Indeed, between these two alternatives, $uv$ versus $u$, it is eventually consistently one or the other, which can be seen as follows:

Write $\overleftarrow{(\omega_n)'_n}=(\xi_n)_n$ and $\overleftarrow{\left(\overline{(\omega_n)_n}\right)'}=(\tau_n)_n$, with $\tau_n=v_{n,1}v_{n,2}\cdots v_{n,m_n}$ for all $n$.\linebreak We claim that there is no index $n$ such that alternative $uv$ at level $n$ is followed by alternative $u$ at level $n+1$.
For if  $\xi_n=v_{n,1}v_{n,2}\cdots v_{n,m_{n-1}}/v_{n,m_n}$ and $\xi_{n+1}=
v_{n+1,1}v_{n+1,2}\cdots v_{n+1,m_{n+1}}/v_{n+1,m_{n+1}+1}$ for some $n$ and some $v_{n+1,m_{n+1}+1}$,   then $\phi_n(\xi_{n+1})=\xi_n$, while $\phi_n(v_{n+1,m_{n+1}})=v_{n,m_n}$ by Remark~\ref{last}. But this is not consistent with Definition~\ref{formal}.

(c) We may assume, without loss of generality (cf.\@ Remark~\ref{classes}), that $(\omega_n)_n$ is of terminating type. Then $\omega_n=v_{n,1}v_{n,2}\cdots  v_{n,{m_n}}$ for all $n\geqslant N$ and either  $\xi_n=v_{n,1}v_{n,2}\cdots v_{n,{m_n-1}}/v_{n,{m_n}}$ for all $n\geqslant N$ or $\xi_n=v_{n,1}v_{n,2}\cdots v_{n,{m_n}}/v_{n,{m_n}+1}$  for all $n\geqslant N$ and some $v_{n,{m_n}+1}$.
Either way, since $\phi_n(v_{n+1,m_{n+1}})=v_{n,{m_n}}$ for all $n\geqslant N$, we have $j_k=m_k$ at every sufficiently large level in Definition~\ref{CompletionDef} for $\overline{(\xi_n)_n}$. Therefore, the value of $drc_n^k$ in Definition~\ref{CompletionDef} is the same for both sequences $(\omega_n)_n$ and $(\xi_n)_n$, so that $\overline{(\omega_n)_n}=\overline{(\xi_n)_n}$.
\end{proof}

\begin{definition}[Words spelled by paths: $\alpha_n\mapsto \omega_n(\alpha_n)\;$] \label{spelling} Given a continuous path \linebreak $\alpha_n:[0,1]\rightarrow X_n$ with $\alpha_n(0)=x_n$, we  let $\omega_n(\alpha_n)\in {\mathcal W}_n$ denote the word ``spelled''  by $\alpha_n$. Specifically, let \[0=s_1\leqslant t_1<s_2\leqslant t_2<\cdots<s_k\leqslant t_k\leqslant1\] be the unique subdivision of $[0,1]$ such that
\begin{itemize}\item[]$\alpha_n(s_i)=\alpha_n(t_i)\in V_n$ for all $1\leqslant i\leqslant k$;\item[] $\alpha_n(u)\cap V_n\subseteq\{\alpha_n(s_i)\}$ for all $1\leqslant i\leqslant k$ and all $u\in [s_i,t_i]$;\item[] $\alpha_n(u)\not\in V_n$ for all $u\not\in \bigcup_{i=1}^k[s_i,t_i]$;  \item[] $\alpha_n(t_i)\not=\alpha_n(s_{i+1})$ for all $1\leqslant i\leqslant k-1$.\end{itemize}
Put $v_i=\alpha_n(s_i)$. If $\alpha_n(1)=v_k$ we define $\omega_n(\alpha_n)=v_1v_2\cdots v_k$, otherwise we put $\omega_n(\alpha_n)=v_1v_2\cdots v_k/u$ where $\alpha_n(1)$ lies on the edge $\{v_k,u\}\in E_n$. \end{definition}

\begin{remark}\label{straight} The word $\omega_n(\alpha_n)$ records the traversed vertices of the edge-path in $X_n$ obtained by straight-line homotopies of $\alpha_n$ on the above subdivision intervals.
\end{remark}

\begin{remark}
  For two paths $\alpha_n,\beta_n:[0,1]\rightarrow X_n$ with $\alpha_n(0)=\alpha_n(1)=\beta_n(0)=x_n$, we have $\omega_n(\alpha_n\cdot \beta_n)=\omega_n(\alpha_n)\omega_n(\beta_n)$.
\end{remark}

Word sequences generated by continuous paths in $X$ are complete:

\begin{lemma}\label{obvious}
For every continuous path  $\alpha=(\alpha_n)_n:([0,1],0)\rightarrow (X,x)$ we have $(\omega_n(\alpha_n))_n\in {\mathcal W}$ and $\overline{(\omega_n(\alpha_n))_n}=(\omega_n(\alpha_n))_n$.
\end{lemma}

\begin{proof}
This follows directly from Definitions~\ref{spelling} and \ref{CompletionDef}, and the continuity of~$\alpha$.
\end{proof}

Conversely, Proposition~\ref{realization} states that every completed word sequence can be realized by a continuous path in $X$. The proof is based on the following lemma.

\begin{lemma}\label{Cauchy}
Given continuous functions $\beta_n:[0,1]\rightarrow X_n$ with the property that $\beta_n$ and $f_n\circ \beta_{n+1}$ are contiguous in $X_n$, the  limits \[ \alpha_n=\lim_{k\rightarrow \infty} f_n\circ f_{n+1}\circ\cdots \circ f_{k-1}\circ \beta_k:[0,1]\rightarrow X_n\] exist and combine to a continuous
    function $\alpha=(\alpha_n)_n:[0,1]\rightarrow X$.
\end{lemma}

\begin{proof}
By Lemma~\ref{setup}, the sequence $(f_n\circ f_{n+1}\circ\cdots \circ f_{k-1}\circ \beta_k:[0,1]\rightarrow X_n)_k$ is uniformly Cauchy.
\end{proof}

\begin{proposition}\label{realization} For every word sequence $(\xi_n)_n\in {\mathcal W}$ or $(\xi_n)_n\in \Omega$, there is a continuous path or loop, respectively, $\alpha=(\alpha_n)_n:([0,1],0)\rightarrow (X,x)$ such that  $(\omega_n(\alpha_n))_n=\overline{(\xi_n)_n}$.
\end{proposition}

\begin{proof} We construct $\alpha$ in the obvious canonical way.
First, we define a piecewise linear continuous path $\beta_1:[0,1]\rightarrow X_1$ based on
 the word $\xi_1=v_1v_2\cdots v_k/\ast$. Let
$0=s_1< t_1<s_2< t_2<\cdots<s_k< t_k=1$ be the partition that subdivides $[0,1]$ into $2k-1$ intervals of equal length and let $\beta_1$ be the unique piecewise linear function on this subdivision with $\beta_1(s_i)=\beta_1(t_i)=v_i$ for all $1\leqslant i \leqslant k$.

 We then define a piecewise linear continuous path $\beta_2:[0,1]\rightarrow X_2$ as follows. Say, $\xi_2=u_1u_2\cdots u_m/\ast$. Let $i_1$ be the maximal index such that $\phi_1(u_1)=\phi_1(u_{i_1})=v_1$ and $\phi_1(u_1u_2\cdots u_{i_1})=v_1$. Subdividing the interval $[s_1,t_1]$ into $2i_1-1$ subintervals of equal length, we define  $\beta_2$ to be alternately constant and linear on these subintervals, the constant values being the vertices $u_1,u_2,\dots, u_{i_1}$. Next,
 let $i_2$ be the maximal index such that $\phi_1(u_{i_1+1}u_{i_1+2}\cdots u_{i_2})$ is the empty word. Subdividing $[t_1,s_2]$ into $2(i_2-i_1)+1$ subintervals of equal length, we define $\beta_2$ to be alternately linear and constant on these subintervals, the constant values being the vertices $u_{i_1+1},u_{i_1+2},\dots, u_{i_2}$. We process the remaining intervals $[s_2,t_2],[t_2,s_3],\dots,[s_k,t_k]$ analogously until $\beta_2$ is fully defined.

Continuing in this fashion, we obtain a sequence $(\beta_n:[0,1]\rightarrow X_n)_n$ of continuous functions  such that $\beta_n$ and $f_n\circ \beta_{n+1}$ are contiguous. Let $\alpha:[0,1]\rightarrow X$ be the limit path provided by Lemma~\ref{Cauchy}. The fact that $\alpha$  has the desired properties follows from the proof of Lemma~\ref{completion}.
\end{proof}

\begin{remark}\label{GroupDiagram}
By Lemma~\ref{obvious} and Proposition~\ref{realization},  $(\alpha_n)_n\mapsto (\omega_n(\alpha_n))_n$ defines a surjection from the set of all continuous loops $L(X,x)$ in $X$, based at $x$, onto the set $\overline{\Omega}$ of all completed word sequences in $\Omega$.
On one hand, the fundamental group $\pi_1(X,x)$ is the image of $L(X,x)$ under the function $(\alpha_n)_n\mapsto [(\alpha_n)_n]$ which forms the homotopy classes. On the other hand, by Lemma~\ref{setequal} and Lemma~\ref{period}(a), we have a surjection from $\overline{\Omega}$ onto the set
$\mathcal G$ of locally eventually constant sequences in $G$ given by $(\omega_n)_n\mapsto (\omega_n)'_n$. In order to circumvent a systematic discussion of the combinatorial relationship between word sequences that represent homotopic paths, we will shift our focus to the function $\varphi:\pi_1(X,x)\rightarrow \mathcal G\subseteq G$ given by $\varphi([(\alpha_n)_n])=(\omega_n(\alpha_n))'_n$, which makes the following diagram commute and which will be shown to be a well-defined isomorphism in Lemma~\ref{cech} and Theorem~\ref{group}.

\hspace{-18pt} \parbox{8in}{
\[
\xymatrix{
L(X,x) \ar[d]_{[\;\cdot\;]}   \ar[rrr]_{surjection}^{(\alpha_n)_n\mapsto(\omega_n(\alpha_n))_n} &&&  \;\overline{\Omega}\;  \ar[d]^{'} & \;\; \overline{\overleftarrow{\mathcal G}} \ar@{_{(}->}[l] \ar@{.>}[dr]_{'}^{\; bijections} &&  \overleftarrow{\mathcal G}\; \ar@{.>}_{``\overline{\;\;\;}"}[ll] \ar@{^{(}->}[r] &\Omega\subseteq  {\mathcal W}\\
\pi_1(X,x) \ar@{.>}[rrr]^{\varphi}_{isomorphism} &&& \; {\mathcal G} \;  & & \; {\mathcal G} \;\;\ar@{^{(}->}[rr] \ar@{.>}_{``\overleftarrow{\;\;\;}"}[ur] & & G\subseteq R}
\]}
So, at the level of word sequences, we obtain the correspondence  $\pi_1(X,x)\cong \overleftarrow{\mathcal G}$.
 In Theorem~\ref{image} (and Remark~\ref{TreeDiagram}), we will establish the more general correspondence between the homotopy classes of paths in $X$ which start at $x$, denoted by $\tilde{X}$, and the elements of the set $\dot{\overleftarrow{\mathcal R}}\subseteq {\mathcal W}$.
 By Lemma~\ref{CurtisFort}, $\tilde{X}$ is a uniquely arcwise connected space. In Corollary~\ref{coarc}, we  show that the radial arcs of $\tilde{X}$, when projected into $X$, precisely spell out the completions of the elements in $\dot{\overleftarrow{\mathcal R}}$. We now work out the details.
\end{remark}

\begin{definition}[Lifts] \label{lifts} Let $\tilde{X}$ denote the set of all homotopy classes $[\alpha]$ of continuous paths $\alpha:([0,1],0)\rightarrow (X,x)$ and let $\tilde{x}$ denote the class containing the constant path. Endow $\tilde{X}$ with the topology generated by the basis comprised of all sets of the form $B([\alpha],U)=\{[\beta]\mid [\beta]=[\alpha\cdot \gamma], \gamma:[0,1]\rightarrow U\}$. Since $X$ is path-connected, we have that  $\tilde{X}$ is path-connected, locally path-connected and metrizable~\cite{FZ2}. Define the map
 $q=(q_n)_n:(\tilde{X},\tilde{x})\rightarrow (X,x)$ by $q([\alpha])=\alpha(1)$, i.e., $q_n([(\alpha_n)_n])=\alpha_n(1)$. Express the elements of the universal covering spaces  $\tilde{X}_n$ of $X_n$ as homotopy classes of continuous paths in $X_n$ starting at $x_n$, i.e., $\tilde{X}_n=\{[\alpha_n]\mid\alpha_n:([0,1],0)\rightarrow (X_n,x_n)\}$,  and let $\tilde{x}_n\in \tilde{X}_n$ denote the class containing the constant path. The covering maps $p_n:\tilde{X}_n\rightarrow X_n$ are given by $p_n([\alpha_n])=\alpha_n(1)$. Lift the given bonding maps $f_n:(X_{n+1},x_{n+1})\rightarrow (X_n,x_n)$ to maps $\tilde{f}_n:(\tilde{X}_{n+1},\tilde{x}_{n+1})\rightarrow (\tilde{X}_n,\tilde{x}_n)$ such that $p_n\circ\tilde{f}_n
=f_n\circ p_{n+1}$ for all $n$. Specifically, $\tilde{f}_n([\alpha_{n+1}])=[f_n\circ \alpha_{n+1}]$. Finally, define $\tilde{q}_n:(\tilde{X},\tilde{x})\rightarrow (\tilde{X}_n,\tilde{x}_n)$ by $\tilde{q}_n([(\alpha_n)_n])=[\alpha_n]$. Then $\tilde{q}_n$ is continuous, $p_n\circ \tilde{q}_n=q_n$ and $\tilde{f}_n\circ \tilde{q}_{n+1}=\tilde{q}_n$ for all~$n$:
\vspace{-12pt}

\[
\xymatrix{  & & & \tilde{X}_{n+1} \ar[dd]^<(.3){p_{n+1}} \ar[dl]^{\tilde{f}_n}\\
\tilde{X}\ar[ddrr]^(.55){q_n} \ar[drrr]_<(.4){q_{n+1}}|<(.56){\hole} \ar[dd]_<(.45){q} \ar[rr]^{\tilde{q}_n } \ar[rrru]^{\tilde{q}_{n+1}}& & \tilde{X}_n \ar[dd]^<(.45){p_n}& \\
& &  & X_{n+1} \ar[dl]^{f_n}\\
X \ar[rrru]^<<<<<<<<{proj}|<(.42){\hole}|<(.56){\hole} \ar[rr]_{proj}& &  X_n & }
\]

\end{definition}

The following fact has essentially been known since \cite{CF1}. We sketch a proof using the argument given there for the Menger cube.

\begin{lemma}\label{CurtisFort} The space $\tilde{X}$ is uniquely arcwise connected and
the map \[\tilde{q}=(\tilde{q}_n)_n:\tilde{X}\rightarrow \hat{X}=\lim_{\longleftarrow} \left(\tilde{X}_1\stackrel{\tilde{f}_1}{\longleftarrow} \tilde{X}_2 \stackrel{\tilde{f}_2}{\longleftarrow} \tilde{X}_3 \stackrel{\tilde{f}_3}{\longleftarrow} \cdots\right)\] sending $[(\alpha_n)_n]\mapsto ([\alpha_n])_n$ is injective.
\end{lemma}

\begin{proof}(Based on \cite{CF1}.)
Since each $\tilde{X}_n$ is a tree, the inverse limit $\hat{X}$ does not contain any simple closed curves. Therefore, every compact path-connected and locally path-connected subspace of $\hat{X}$ is a dendrite and hence contractible. Consequently, the map $\displaystyle\tilde{q}=(\tilde{q}_n)_n:\tilde{X}\rightarrow \hat{X}$ is injective and $\tilde{X}$ contains no simple closed curve.
\end{proof}

\begin{remark}\label{simply}
The map  $\displaystyle\tilde{q}=(\tilde{q}_n)_n:\tilde{X}\rightarrow \hat{X}$ is always well-defined and continuous for any inverse limit $X$ of topological spaces $X_n$, even if $\tilde{X}$ is not simply connected. (This follows directly from the definition of the topologies on $\tilde{X}$ and $\tilde{X}_n$.) However,  if $\tilde{q}$ happens to be injective and if each $X_n$ is a compact metric space, then the natural homomorphism   $\tilde{q}|_{\pi_1(X,x)}:\pi_1(X,x)\rightarrow \check{\pi}_1(X,x)$ into the first \v{C}ech \linebreak homotopy group $\displaystyle \check{\pi}_1(X,x)=\lim_{\longleftarrow} \left( \pi_1(X_1,x_1) \stackrel{f_{1\#}}{\leftarrow}\pi_1(X_2,x_2) \stackrel{f_{2\#}}{\leftarrow}\pi_1(X_3,x_3) \stackrel{f_{3\#}}{\leftarrow}\cdots\right)$  is injective so that $\tilde{X}$ is simply connected \cite{FZ2}.
\end{remark}

\begin{definition} \label{phidef}
We define functions $\varphi_n:\tilde{X}_n\rightarrow {\mathcal W}_n'$ by $\varphi_n([\alpha_n])=\omega_n(\alpha_n)'$.
\end{definition}

We record the following straightforward lemma without proof:

  \begin{lemma} \label{cech} Each $\Omega_n'$ forms a free group under the operation $\omega_n\ast\xi_n=(\omega_n\xi_n)'$ and $\varphi_n:\pi_1(X_n,x_n)\rightarrow \Omega_n'$  is an isomorphism. Moreover, the following diagrams commute for all $n$:
\[\xymatrix{ & \pi_1(X,x) \ar[ld]_{\tilde{q}_n} \ar[rd]^{\tilde{q}_{n+1}} & \\ \pi_1(X_n,x_n)\ar[d]_{\varphi_n} & & \pi_1(X_{n+1},x_{n+1}) \ar[ll]_{f_{n\#}} \ar[d]^{\varphi_{n+1} }\\
\Omega_n' &  & \Omega_{n+1}'\ar[ll]_{\phi'_n } }\]
Consequently,  $\check{\pi}_1(X,x)\cong G$ are isomorphic and the function $\varphi:\pi_1(X,x)\rightarrow  G$, given by $\varphi([(\alpha_n)_n])=(\varphi_n([\alpha_n]))_n=(\omega_n(\alpha_n))'_n$,
 is an injective homomorphism.
\end{lemma}

\begin{theorem}[Theorem~A]\label{group}
 We have $\varphi(\pi_1(X,x))={\mathcal G}\subseteq G$. Hence,  $\overleftarrow{\mathcal G}$ forms a group under the operation $\ast$, given by $(w_n)_n \ast (\xi_n)_n=\overleftarrow{(\omega_n\xi_n)'_n}$, and  $\overleftarrow{\mathcal G}\cong\pi_1(X,x)$.
\end{theorem}

\begin{proof}
Let $L(X,x)$ denote the set of all continuous loops in $X$ which are based at~$x$. For a given $(\alpha_n)_n\in L(X,x)$, we have $\varphi([(\alpha_n)_n])=(\omega_n(\alpha_n))_n'\in \Omega'={\mathcal G}$ by Lemma~\ref{setequal}.
Conversely, let $(g_n)_n\in {\mathcal G}$ be given.
By Proposition~\ref{realization}, there is an $(\alpha_n)_n\in L(X,x)$ with $(\omega_n(\alpha_n))_n=\overline{\overleftarrow{(g_n)_n}}$. Then, by Lemma~\ref{period}(a) and Remark~\ref{representation}, we have \[\varphi([(\alpha_n)_n])=(\omega_n(\alpha_n))_n'=
\left(\overline{\overleftarrow{(g_n)_n}}\right)'=\left(\overleftarrow{(g_n)_n}\right)'=(g_n)_n.\]
Therefore, $\varphi(\pi_1(X,x))={\mathcal G}$. The fact that the operation $\ast$ corresponds to multiplication in $\pi_1(X,x)$ is verified in more generality in Remark~\ref{action} below.

\end{proof}

\begin{lemma}\label{phi}
The following diagrams commute for all n:
\[ \xymatrix{& \tilde{X}\ar[ld]_{\tilde{q}_n} \ar[rd]^{\tilde{q}_{n+1}} & \\ \tilde{X}_n \ar[d]_{\varphi_n}& &\tilde{X}_{n+1} \ar[ll]_{\tilde{f}_n} \ar[d]^{\varphi_{n+1} }\\
{\mathcal W}_n' & & {\mathcal W}_{n+1}'\ar[ll]_{\phi'_n } }\]
The combined functions $\displaystyle (\varphi_n)_n:\hat{X}\rightarrow  R$
    yield an injection so that we obtain an injective function  $\varphi:\tilde{X}\hookrightarrow  R$ given by $\varphi([(\alpha_n)_n])=(\varphi_n([\alpha_n]))_n=(\omega_n(\alpha_n))'_n$.
\end{lemma}

\begin{proof}
The difference between this lemma and Lemma~\ref{cech} is the fact that the functions $\varphi_n:\tilde{X}_n\rightarrow {\mathcal W}_n'$ are not bijective.  For a reduced word $r_n\in{\mathcal W}_n'$ of the form $r_n=v_1v_2\cdots v_k$, the preimage $\varphi_n^{-1}(\{r_n\})$ is a vertex of the tree $\tilde{X}_n$. For  a reduced word  $r_n\in{\mathcal W}_n'$ of the form  $r_n=v_1v_2\cdots v_k/v_{k+1}$, the preimage $\varphi_n^{-1}(\{r_n\})$  is a half-open edge of the tree $\tilde{X}_n$, where $v_1v_2\cdots v_k$ is the shortest edge-path representative for the homotopy class of paths that connect $x_n$ to the included vertex of the projection of this edge in $X_n$. Therefore, by Lemma~\ref{setup}, the combined functions  $(\varphi_n)_n$ yield an  injection. Note that although each $\varphi_n$ is clearly surjective, the combined functions $\displaystyle (\varphi_n)_n:\hat{X}\rightarrow  R$ need not yield a surjection, as illustrated in the following remark.
\end{proof}

\begin{remark}\label{nonsurj} An example for which $\displaystyle (\varphi_n)_n:\hat{X}\rightarrow  R$ is not surjective is given by $X=[0,1]$, expressed as an inverse limit of subdivisions $X_n$ of $[0,1]$ with $f_n=id|_{X_n}$ and $x_n=0$ for all $n$. Label the vertices  of $X_n$  as $0=v_{n,1}<v_{n,2}<\cdots< v_{n,m_n}=1$, and form the word $r_n=v_{n,1}v_{n,2}\cdots v_{n,{m_n-1}}/v_{m_n}\in {\mathcal W}'_n$. Then $(r_n)_n\in R$. However, $(r_n)_n$  is not in the image of $\displaystyle (\varphi_n)_n:\hat{X}\rightarrow  R$.
\end{remark}

\begin{theorem}\label{image}
We have\vspace{5pt}

 \begin{itemize}
\item[(a)] $\varphi(\tilde{X})\subseteq {\mathcal R}\subseteq R$;\vspace{10pt}

 \item[(b)] $\overleftarrow{\varphi(\tilde{X})}=\dot{\overleftarrow{\mathcal R}}\subseteq {\mathcal W}$, yielding a bijective correspondence between $\tilde{X}$ and  $\dot{\overleftarrow{\mathcal R}}$;

 \item[(c)] $\left(\overline{\overleftarrow{\varphi(\tilde{y})}}\right)'=\varphi(\tilde{y})$ for all $\tilde{y}\in \tilde{X}$.
 \end{itemize}
\end{theorem}

\begin{proof}
The proof is similar to that of Theorem~\ref{group}.

(a) This follows from Lemma~\ref{setequal}.

(b) Put ${\mathcal S}=\overleftarrow{\varphi(\tilde{X})}$. While not every element of $\mathcal S$ is of terminating type, it follows from Definition~\ref{comp} and Definition~\ref{spelling}  that $\dot{\mathcal S}={\mathcal S}$. Hence, $\overleftarrow{\varphi(\tilde{X})}\subseteq \dot{\overleftarrow{\mathcal R}}$ by Part~(a).  For the reverse inclusion, let $(r_n)_n\in {\mathcal R}$.
By  Proposition~\ref{realization} we may choose $\tilde{y}=[(\alpha_n)_n]\in \tilde{X}$ with
$(\omega_n(\alpha_n))_n=\overline{\overleftarrow{(r_n)_n}}$. Then $\overleftarrow{\varphi(\tilde{y})}\doteq \overleftarrow{(r_n)_n}$, as in the proof of Theorem~\ref{group}, but using  Lemma~\ref{period}(b) instead of Lemma~\ref{period}(a):
\[\overleftarrow{\varphi(\tilde{y})}=\overleftarrow{\varphi([(\alpha_n)_n])}=\overleftarrow{(\omega_n(\alpha_n))_n'}=
\overleftarrow{\left(\overline{\overleftarrow{(r_n)_n}}\right)'}\doteq \overleftarrow{\left(\overleftarrow{(r_n)_n}\right)'}=\overleftarrow{(r_n)_n}.\]
The fact that $\tilde{z}\mapsto \overleftarrow{\varphi(\tilde{z})}$ defines a bijection from $\tilde{X}$ onto $\dot{\overleftarrow{\mathcal R}}$ follows now from Lemma~\ref{phi} and Remark~\ref{representation}.

(c)
Let $\tilde{y}\in \tilde{X}$. Then, as in the proof of Part~(b), we have $\overleftarrow{\left(\overline{\overleftarrow{\varphi(\tilde{y})}}\right)'}\doteq \overleftarrow{\varphi(\tilde{y})}$. Also, by Proposition~\ref{realization}, $\left(\overline{\overleftarrow{\varphi(\tilde{y})}}\right)'\in \varphi(\tilde{X})$. Hence, $\left(\overline{\overleftarrow{\varphi(\tilde{y})}}\right)'=\varphi(\tilde{y})$ by Part~(b).
\end{proof}

\begin{remark}\label{TreeDiagram} By Theorem~\ref{image} and Lemma~\ref{period}(c), we now have the following commutative diagram, where $P(X,x)$ denotes all continuous paths in $X$  which start at $x$.

\hspace{-18pt} \parbox{8in}{
\[
\xymatrix{
P(X,x) \ar[d]_{[\;\cdot\;]}   \ar[rrr]_{surjection}^{(\alpha_n)_n\mapsto(\omega_n(\alpha_n))_n} &&&  \;\overline{\mathcal W}\;  \ar[d]^{'} & \;\; \overline{\overleftarrow{\mathcal R}} \ar@{_{(}->}[l] \ar@{.>}[dr]_{'}^{\;\;\;\;\; bijections} &&  \dot{\overleftarrow{\mathcal R}}\; \ar@{.>}_{``\overline{\;\;\;}"}[ll] \ar@{^{(}->}[r] & {\mathcal W}\\
\;\tilde{X}\; \ar@{>->}[rrr]^{\varphi}_{injection} &&& \; {\mathcal R} \;  & & \; \varphi(\tilde{X}) \;\;\ar@{^{(}->}[rr] \ar@{.>}_{``\overleftarrow{\;\;\;}"}[ur] & & R}
\]}
\end{remark}

\begin{remark} \label{action} Under the correspondence of Theorem~\ref{image}(b) between $\tilde{X}$ and $\dot{\overleftarrow{\mathcal R}}$, the natural action of $\pi_1(X,x)$ on $\tilde{X}$, given  by $[\alpha].[\beta]=[\alpha\cdot\beta]$, corresponds to the action of $\overleftarrow{\mathcal G}$ on $\dot{\overleftarrow{\mathcal R}}$, given by $(\omega_n)_n.(\xi_n)_n=\overleftarrow{(\omega_n\xi_n)_n'}$. To see this, suppose $(\omega_n)_n=\overleftarrow{\varphi(\tilde{y})}$ and $(\xi_n)_n=\overleftarrow{\varphi(\tilde{z})}$, with $\tilde{y}=[(\alpha_n)_n]$ and $\tilde{z}=[(\beta_n)_n]$. Then  we have $(\omega_n)_n=\overleftarrow{(\omega_n(\alpha_n))'_n}$ so that $\omega_n'=\omega_n(\alpha_n)'$ for all $n$.
Similarly, $\xi_n'=\omega_n(\beta_n)'$ for all $n$. Hence, $\overleftarrow{(\omega_n\xi_n)'_n}=\overleftarrow{(\omega_n'\xi_n')'_n}=\overleftarrow{(\omega_n(\alpha_n)'\omega_n(\beta_n)')'_n}=
\overleftarrow{(\omega_n(\alpha_n)\omega_n(\beta_n))'_n}=\overleftarrow{(\omega_n(\alpha_n\cdot \beta_n))'_n}=\overleftarrow{\varphi([(\alpha_n\cdot\beta_n)_n])}
=\overleftarrow{\varphi([(\alpha_n)_n].[(\beta_n)_n])}
=\overleftarrow{\varphi(\tilde{y}.\tilde{z})}$.
\end{remark}

\begin{corollary}[Theorem~D]\label{ThmD}
The group $\overleftarrow{\mathcal G}\cong \pi_1(X,x)$ acts freely and by homeomorphism on the $\mathbb{R}$-tree $\dot{\overleftarrow{\mathcal R}}$ via its natural action $(\omega_n)_n.(\xi_n)_n=\overleftarrow{(\omega_n \xi_n)'_n}$.
\end{corollary}

\begin{proof} Based on Remark~\ref{simply}, we may apply \cite[Theorem~4.10]{FZ2}. In particular, $\pi_1(X,x)$ acts freely and by homeomorphism on the generalized universal covering space $\tilde{X}$. The results now follow from Theorem~\ref{image}(b) and Remark~\ref{action}.\end{proof}

\begin{definition} For $\tilde{y},\tilde{z}\in \tilde{X}$, we denote   the unique arc in $\tilde{X}$ from $\tilde{y}$ to $\tilde{z}$ by $[\tilde{y},\tilde{z}]$.
\end{definition}

Corollary~\ref{coarc} below states that the arc $[\tilde{x},\tilde{y}]$ in $\tilde{X}$ from the base point $\tilde{x}$  to a point $\tilde{y}$, when projected into the approximating graphs $X_n$ of $X$, spells out the word sequence $\overline{\overleftarrow{\varphi(\tilde{y})}}$. The proof follows from Proposition~\ref{arc}, which in turn uses the following:

\begin{lemma}\label{sub} Let $(\omega_n)_n,(\xi_n)_n\in \mathcal W$. Express the word sequences $\overline{(\omega_n)_n}, \overline{(\xi_n)_n}\in {\mathcal W}$ as $\overline{(\omega_n)_n}=(\bar{\omega}_n)_n$ and
$\overline{(\xi_n)_n}=(\bar{\xi}_n)_n$ with $\bar{\omega}_n, \bar{\xi}_n\in {\mathcal W}_n$. If $\omega_n$ is a subsequence of $\xi_n$ for all $n$, then $\bar{\omega}_n$ is a subsequence of $\bar{\xi}_n$ for all $n$.
\end{lemma}

\begin{remark} The notion of subsequence is in the sense of the proof of Lemma~\ref{setequal}.
\end{remark}

\begin{proof}  $drc_n^k$ for $(\omega_n)_n$ produces subsequences of $drc_n^k$ for $(\xi_n)_n$.
\end{proof}

\begin{proposition}\label{arc}
Let $\tilde{y}\in \tilde{X}$, let  $\tilde{\beta}:[0,1]\rightarrow \tilde{X}$ be a parametrization of the arc $[\tilde{x},\tilde{y}]$ in $\tilde{X}$ and put $(\beta_n)_n=q\circ\tilde{\beta}:[0,1]\rightarrow X$. Then
\[\overline{\overleftarrow{(\omega_n(\beta_n))'_n}}=(\omega_n(\beta_n))_n.\]
\end{proposition}

\begin{proof} Fix $k\geqslant 1$. Note that  $\tilde{y}=\tilde{\beta}(1)=[q\circ \tilde{\beta}]=[(\beta_n)_n]$.
Put $(\xi_n)_n=\overleftarrow{(\omega_n(\beta_n))'_n}$ and
 express $\overline{(\xi_n)_n}=(\bar{\xi}_n)_n$.
  Since $(\xi_n)_n$ is the stabilization of  $\varphi(\tilde{y})=(\omega_n(\beta_n))'_n$, we have $\xi_k=\phi_k\circ\phi_{k+1}\circ\cdots\circ\phi_{n-1}(\omega_n(\beta_n)')$, for all sufficiently large $n$, which is a subsequence of
$\phi_k\circ\phi_{k+1}\circ\cdots\circ\phi_{n-1}(\omega_n(\beta_n))=\omega_k(\beta_k)$. By Lemma~\ref{obvious}, we have $\overline{(\omega_n(\beta_n))_n}=(\omega_n(\beta_n))_n$, so that $\bar{\xi}_k$ is a subsequence of $\omega_k(\beta_k)$ by Lemma~\ref{sub}.

By Proposition~\ref{realization}, there is a path $\gamma=(\gamma_n)_n:([0,1],0)\rightarrow (X,x)$ such that $(\omega_n(\gamma_n))_n=\overline{(\xi_n)_n}$. Let $\tilde{\gamma}:([0,1],0)\rightarrow (\tilde{X},\tilde{x})$ be the lift with $q\circ \tilde{\gamma}=\gamma$. Then, by Lemma~\ref{period}(b),

\[\overleftarrow{\varphi(\tilde{\gamma}(1))}=\overleftarrow{\varphi([\gamma])}=\overleftarrow{\varphi([(\gamma_n)_n])}=\overleftarrow{(\varphi_n([\gamma_n]))_n}
=\overleftarrow{(\omega_n(\gamma_n)_n)'_n}=\] \[=\overleftarrow{\left(\overline{(\xi_n)_n}\right)'}\doteq \overleftarrow{(\xi_n)'_n}=\overleftarrow{\left(\overleftarrow{(\omega_n(\beta_n))'_n}\right)'}=\overleftarrow{(\omega_n(\beta_n))'_n}=
\overleftarrow{\varphi(\tilde{y})}.\]

Therefore, by Lemma~\ref{period}(c) and Theorem~\ref{image}(c), $\varphi(\tilde{\gamma}(1))=\varphi(\tilde{y})$. Hence, by Lemma~\ref{phi}, the path
$\tilde{\gamma}:[0,1]\rightarrow \tilde{X}$ connects the endpoints of the arc $\tilde{\beta}([0,1])$. Since $\tilde{X}$ is uniquely arcwise connected, this implies (directly from the definition) that $\omega_k(\beta_k)$ is a subsequence of $\omega_k(\gamma_k)=\bar{\xi}_k$.

Hence,  $\omega_k(\beta_k)=\bar{\xi}_k$, each being a subsequence of the other.
\end{proof}

\begin{corollary}\label{coarc}
Let $\tilde{y}=[(\alpha_n)_n)]\in \tilde{X}$, let $\tilde{\beta}:[0,1]\rightarrow \tilde{X}$ be a parametrization of the arc $[\tilde{x},\tilde{y}]$ in $\tilde{X}$ and put $(\beta_n)_n=q\circ \tilde{\beta}:[0,1]\rightarrow X$. Then
\[\overline{\overleftarrow{\varphi(\tilde{y})}}=\overline{\overleftarrow{(\omega_n(\alpha_n))'_n}}=(\omega_n(\beta_n))_n.\]
\end{corollary}

\begin{proof} As in the proof of Proposition~\ref{arc}, we have $[(\alpha_n)_n]=\tilde{y}=[(\beta_n)_n]$. Hence,
$(\omega_n(\alpha_n))'_n=
\varphi(\tilde{y})=(\omega_n(\beta_n))'_n$.
\end{proof}

\begin{example}\label{L-space} Note that the stabilization of $(\omega_n(\alpha_n))'_n$ in Corollary~\ref{coarc} might not be complete. This can be observed, for example, in the one-point compactification $\mathbb{L}$ of an infinitely long ladder, expressed as the limit of an inverse sequence satisfying Lemma~\ref{setup}. Such a space and its defining sequence are depicted in Figure~\ref{L-picture}.
Let $x=a$ be the base point of the space $X=\mathbb{L}$  and let  $\tilde{y}$ be the homotopy class of a parametrization $(\alpha_n)_n$ of the arc $ab\cup bc$. While the words of the sequence $\overleftarrow{(\omega_n(\alpha_n))'_n}=(\omega_n(\alpha_n))'_n$ never include the top vertex of the corresponding approximating graph, all words of the completion do. \qed
\begin{figure}[h]
\fbox{\parbox{4.5in}{\hspace{.25in}
\includegraphics{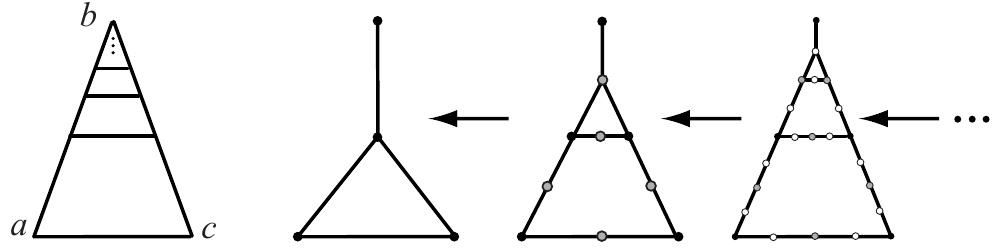}}}
\caption{\label{L-picture} The space $\mathbb{L}$ (left) and its defining sequence (right). Although all paths in the homotopy class of the arc $ab\cup bc$ travel through the point $b$, the reduced projections of this arc into the approximating spaces yield an already stabilized word sequence, none of whose words include the top vertex. The completion process is designed to remedy such omissions by reinserting these vertices.}
\end{figure}

\end{example}

\begin{corollary}[Theorem~C]\label{differences}
  Let $(\omega_n)_n,(\xi_n)_n \in \overleftarrow{\mathcal G}$. Say, $(\omega_n)_n=\overleftarrow{\varphi(\tilde{y})}$ and $(\xi_n)_n=\overleftarrow{\varphi(\tilde{z})}$ with $\tilde{y},\tilde{z}\in \pi_1(X,x)$. Let $\tilde{\gamma}:[0,1]\rightarrow \tilde{X}$ be a parametrization of the arc $[\tilde{y},\tilde{z}]$ in $\tilde{X}$ and put $(\gamma_n)_n=q\circ \tilde{\gamma}:[0,1]\rightarrow X$. Then $\overline{(\omega_n)_n^{-1}\ast (\xi_n)_n}=(\omega_n(\gamma_n))_n$.
\end{corollary}

\begin{proof} Say, $\tilde{y}=[(\alpha_n)_n]$ and $\tilde{z}=[(\beta_n)_n]$. Put $\alpha^-_n(t)=\alpha_n(1-t)$ and
consider the homeomorphism $\psi:\tilde{X}\rightarrow \tilde{X}$ given by $\psi([(\tau_n)_n])=[(\alpha^-_n\cdot \tau_n)_n]$ (cf.\@ \cite[Lemma~2.6]{FZ2}). Then $\psi\circ \tilde{\gamma}:[0,1]\rightarrow \tilde{X}$  parametrizes the  arc from $\tilde{x}$ to $[(\delta_n)_n]=[(\alpha_n)_n]^{-1}.[(\beta_n)_n]$  and $q\circ\psi\circ \tilde{\gamma}=q\circ \tilde{\gamma}$. On one hand, $\overline{\overleftarrow{\varphi([(\delta_n)_n]}}=(\omega_n(\gamma_n))_n$ by Corollary~\ref{coarc}. On the other hand, $\overleftarrow{\varphi([(\delta_n)_n]}=(\omega_n)_n^{-1}\ast (\xi_n)_n$ by Remark~\ref{action}.
\end{proof}

\begin{definition}
 Given $\tilde{y},\tilde{z}\in\tilde{X}$, we define $\tilde{y}\wedge \tilde{z}\in\tilde{X}$ by $[\tilde{x},\tilde{y}]\cap[\tilde{x},\tilde{z}]=[\tilde{x},\tilde{y}\wedge\tilde{z}]$.
\end{definition}

\begin{corollary}\label{meet} Let $\tilde{y},\tilde{z}\in \tilde{X}$. Then the following hold: \vspace{5pt} \begin{itemize}
\item[(a)] $\overline{\overleftarrow{\varphi(\tilde{y}\wedge\tilde{z})}}\doteq\overline{\overleftarrow{\varphi(\tilde{y})}}\Cap  \overline{\overleftarrow{\varphi(\tilde{z})}}$\vspace{10pt}
    \item[(b)] $\|\overline{\overleftarrow{\varphi(\tilde{y}\wedge\tilde{z})}}\|=\|\overline{\overleftarrow{\varphi(\tilde{y})}}\Cap \overline{\overleftarrow{\varphi(\tilde{z})}}\|$
        \end{itemize}
\end{corollary}

\begin{proof}
Choose any parametrizations $\tilde{\beta}:[0,1]\rightarrow \tilde{X}$ of $[\tilde{x},\tilde{y}]$,  $\tilde{\gamma}:[0,1]\rightarrow \tilde{X}$  of $[\tilde{x},\tilde{z}]$, and  $\tilde{\alpha}:[0,1]\rightarrow \tilde{X}$ of $[\tilde{x},\tilde{y}\wedge\tilde{z}]$.
 Put  $\beta=q\circ \tilde{\beta}=(\beta_n)_n$, $\gamma=q\circ \tilde{\gamma}=(\gamma_n)_n$ and $\alpha=q\circ \tilde{\alpha}=(\alpha_n)_n$,
  so that $\tilde{y}=[(\beta_n)_n]$,  $\tilde{z}=[(\gamma_n)_n]$, $\tilde{y}\wedge \tilde{z}=[(\alpha_n)_n]$. By Corollary~\ref{coarc}, we have $\overline{\overleftarrow{\varphi(\tilde{y}\wedge\tilde{z})}}=(\omega_n(\alpha_n))_n$ and $\overline{\overleftarrow{\varphi(\tilde{y})}}\Cap \overline{\overleftarrow{\varphi(\tilde{z})}}=(\omega_n(\beta_n))_n\Cap (\omega_n(\gamma_n))_n$.
 Let $\tau_n=\phi_n\circ\phi_{n+1}\circ \cdots\circ \phi_{k-1}(\omega_{k}(\beta_k)\cap \omega_{k}(\gamma_k))$ for sufficiently large~$k$. Clearly, $\omega_k(\alpha_k)$ is  an initial substring of both $\omega_k(\beta_k)$ and $\omega_k(\gamma_k)$, and hence of $\omega_k(\beta_k)\cap\omega_k(\gamma_k)$, for all~$k$. Now, $\omega_n(\alpha_n)=\phi_n\circ\phi_{n+1}\circ \cdots\circ \phi_{k-1}(\omega_k(\alpha_k))$, so that $\omega_n(\alpha_n)$ is an initial substring of $\tau_n$ or all $n$. Now fix $n$.
Say, $\omega_n(\alpha_n)=v_1v_2\cdots v_i/\ast$ and $\tau_n=v_1v_2\cdots v_iv_{i+1}\cdots v_j/\ast$. We claim that $i=j$. Suppose, to the contrary, that $i<j$.
Choose $t^\beta$ and $t^\gamma$ such that $\tilde{y}\wedge \tilde{z}=\tilde{\beta}(t^{\beta})=\tilde{\gamma}(t^\gamma)$.
Let $s^{\beta_n}_{i+1}$ and $s^{\gamma_n}_{i+1}$ be the subdivision points of $[0,1]$ which, as in Definition~\ref{spelling}, put the letter $v_{i+1}$ into the words $\omega_n(\beta_n)$ and $\omega_n(\gamma_n)$, respectively. Then $t^\beta<s^{\beta_n}_{i+1}$ and $t^\gamma<s^{\gamma_n}_{i+1}$.
 By definition of $\tau_n$, we must have $\omega_k(\beta_k|_{[0,s^{\beta_n}_{i+1}]})=\omega_k(\gamma_k|_{[0,s^{\gamma_n}_{i+1}]})\in {\mathcal P}_k$ for all $k>n$. Hence, $\tilde{\beta}(s^{\beta_n}_{i+1})=[\beta_k|_{[0,s^{\beta_n}_{i+1}]}]=[\gamma_k|_{[0,s^{\gamma_n}_{i+1}]}]=\tilde{\gamma}(s^{\gamma_n}_{i+1})$ by Lemma~\ref{phi}, which contradicts the choice of $t^\beta$ and $t^\gamma$.
\end{proof}

\begin{example} It is possible that  $\omega_n(\alpha_n)=v_1v_2\cdots v_i$ and $\tau_n=v_1v_2\cdots v_i/v_{i+1}$ in the proof of Corollary~\ref{meet}.
 For example, consider again the space $X=\mathbb{L}$ depicted in Figure~\ref{L-picture}, but this time with base point $x=b$, and let $\tilde{y}$ and $\tilde{z}$ be the homotopy classes of the arcs $ba$ and $bc$, respectively. \qed
\end{example}

The proofs of the following three lemmas are a combination of straightforward inductive arguments, which can be extracted from \cite{MO}, and Corollary~\ref{coarc}. We include them for completeness.

\begin{lemma}\label{one}
Let $(\omega_n)_n\in {\mathcal W}$. If $\omega_n=v_1v_2\cdots v_k/\ast$ and $2\leqslant i\leqslant k$, then \[|v_1|=\frac{1}{2^n}  \mbox{ and } 0<\frac{|v_{i-1}|}{2^n}\leqslant |v_i|\leqslant\frac{1}{2} \left( \frac{3}{4}\right)^{n-1}.\]
\end{lemma}

\begin{proof}
It follows from the recursive step of Definition~\ref{weights} that $|v_1|=\frac{1}{2^n}$ and that all weights are positive rational numbers. We first show, by induction on $n$, that $|v_i|\leqslant\frac{1}{2} \left( \frac{3}{4}\right)^{n-1}$. For $n=1$, we have $|v_i|=\frac{1}{2^i}$. This establishes the base case. Write $\omega_{n+1}=u_1u_2\cdots u_m/\ast$ and inductively assume that $|v_i|\leqslant\frac{1}{2} \left( \frac{3}{4}\right)^{n-1}$ for all $i$. By Definition~\ref{weights}, for each $j\in\{1,2,\ldots,m\}$ there are two cases: either $|u_j|=|v_i|/2^s$ for some $1\leqslant i\leqslant k$ and $s\geqslant 1$, or $|u_j|=|v_{i-1}|/2^s+|v_i|/2$ for some $2\leqslant i \leqslant k$ and $s\geqslant 2$. In the first case, we have $|u_j|\leqslant \frac{1}{2} \left( \frac{3}{4}\right)^{n-1}2^{-s}\leqslant \frac{1}{4} \left( \frac{3}{4}\right)^{n-1}<\frac{1}{2} \left( \frac{3}{4}\right)^{n}$. In the second case, we have $|u_j|\leqslant |v_{i-1}|/4+|v_i|/2\leqslant \left(\frac{1}{4}+\frac{1}{2}\right)\frac{1}{2} \left( \frac{3}{4}\right)^{n-1}= \frac{1}{2} \left( \frac{3}{4}\right)^{n}$.

Finally, we show, by induction on $n$, that $|v_{i-1}|\leqslant 2^n|v_i|$ for all $2\leqslant i \leqslant k$. When $n=1$, we have $|v_{i-1}|=\frac{1}{2^{i-1}}$ and $|v_i|=\frac{1}{2^i}$, establishing the base case. Again, write $\omega_{n+1}=u_1u_2\cdots u_m/\ast$ and inductively assume that $|v_{i-1}|\leqslant 2^n|v_i|$ for all $2\leqslant i \leqslant k$. If $|u_j|=|v_{i-1}|/2^s+|v_i|/2$, then $|u_{j-1}|=|v_{i-1}|/2^s$ so that $|u_j|>|u_{j-1}|$ and $|u_{j-1}|<2^{n+1}|u_j|$. So, we may assume that $|u_j|=|v_i|/2^s$ for some $1\leqslant i\leqslant k$ and $s\geqslant 1$.
Either $|u_{j-1}|=2|u_j|$, or $s=2$ and $|u_{j-1}|=|v_{i-1}|/2^t+|v_i|/2$ for some $t\geqslant 2$. Assuming the second case, as we may, we get from the induction hypothesis that
$|u_{j-1}|\leqslant 2^n|v_i|/2^t+|v_i|/2\leqslant 2^n|v_i|/4+2|v_i|/4=(2^n+2)|u_j|\leqslant 2^{n+1}|u_j|$.
\end{proof}

\begin{lemma}\label{two}
 Let $(\omega_n)_n\in {\mathcal W}$. If $\omega_n=v_1v_2\cdots v_k/\ast$ and $k\geqslant 2$, then  \[|\omega_n|-|v_{k-1}|-|v_k|\leqslant  \|(\omega_m)_m\| < |\omega_n|.\]
\end{lemma}

\begin{proof}
 It suffices to show that $|\omega_n|-|v_{k-1}|-|v_k|<|\omega_m|$ for all $m\geqslant n$. If $m=n$, this is trivial. If $m=n+1$, the inequality holds, because $|\omega_{n+1}|=|\omega_n|-|v_k|/2^i$ for some $i\geqslant 1$, by Definition~\ref{weights}. Now, suppose $m=n+2$ and express $\omega_{n+1}$ as  $\omega_{n+1}=u_1u_2\cdots u_s/\ast$. Then $|\omega_{n+2}|=|\omega_{n+1}|-|u_s|/2^j$ for some $j\geqslant 1$. In turn, either $|u_s|=|v_k|/2^t$ or $|u_s|=|v_{k-1}|/2^t+|v_k|/2$ for some $t\geqslant 2$. So, either $|\omega_{n+2}|=|\omega_n|-|v_k|/2^i-|v_k|/2^{t+j}$ or
$|\omega_{n+2}|=|\omega_n|-|v_k|/2^i-|v_{k-1}|/2^{t+j}-|v_k|/2^{1+j}$. Either way, $|\omega_n|-|v_{k-1}|-|v_k|<|\omega_{n+2}|$.

To see the general induction logic, express $\omega_{m-1}$ as
$\omega_{m-1}=u_1u_2\cdots u_s/\ast$. Then $|\omega_{m}|=|\omega_{m-1}|-|u_s|/2^j$ for some $j\geqslant 1$. Observe that every substring
 in the recursion step of Definition~\ref{weights}, except possibly the last one, has length at least two. Hence, each
 $|u_i|$ is a rational linear combination of $|v_1|, |v_2|, \ldots,|v_k|$ with coefficients in $[0,1)$, such that at most two coefficients are positive and the positive coefficients are consecutive. Moreover,
 the $|v_{k}|$-coefficient for $|u_s|$ must be positive.
\end{proof}

\begin{lemma}\label{three} Let $\tilde{y},\tilde{z}\in \tilde{X}$ such that $\tilde{y}\in [\tilde{x},\tilde{z}]$ and $\tilde{y}\not=\tilde{z}$.
 Say, $\overline{\overleftarrow{\varphi(\tilde{y})}}=(\overline{\omega}_n)_n$ and $\overline{\overleftarrow{\varphi(\tilde{z})}}=(\overline{\xi}_n)_n$. Suppose $n$ is sufficiently large so that
   $\overline{\omega}_n=v_1v_2\cdots v_k/\ast$ and \linebreak $\overline{\xi}_n=v_1v_2\cdots v_k u_1 u_2 \cdots u_m/\ast$ with $m\geqslant 3$. Then the weights of the letters\linebreak  $v_1, v_2, \cdots, v_k$ agree in both word sequences. Moreover,  for every $i\in\{1, 2,\dots, m-2\}$ we have  \[\Big\|\overline{\overleftarrow{\varphi(\tilde{z})}}\Big\|-\Big\|\overline{\overleftarrow{\varphi(\tilde{y})}}\Big\|\geqslant |u_{i}|>0.\]
\end{lemma}

\begin{proof} The fact that $n$ can be chosen as claimed and that the weights of the letters $v_1, v_2, \ldots, v_k$ agree in both word sequences follows from Corollary~\ref{coarc}. Now, fix any $i\in\{1, 2,\dots, m-2\}$. By Lemma~\ref{two}, we have \[\Big\|\overline{\overleftarrow{\varphi(\tilde{z})}}\Big\|\geqslant |\overline{\xi}_n|-|u_{m-1}|-|u_m|\geqslant |v_1|+|v_2|+\cdots +|v_k|+|u_i|= |\overline{\omega}_n|+|u_i|\geqslant \Big\|\overline{\overleftarrow{\varphi(\tilde{y})}}\Big\|+|u_i|.\]
\end{proof}

\begin{remark}\label{increasing}
Clearly, $\Big\|\overline{\overleftarrow{\varphi(\tilde{x})}}\Big\|=0$. Moreover, it follows from Lemma~\ref{three} that for every  $\tilde{z}\in\tilde{X}$, the function $\tilde{y}\mapsto \Big\|\overline{\overleftarrow{\varphi(\tilde{y})}}\Big\|$ is increasing on the arc $[\tilde{x},\tilde{z}]$.
\end{remark}

\begin{remark}\label{ex}
Let $S$ be a set with base point $s_0\in S$ and let $\tau$ be a topology on $S$ such that the space $(S,\tau)$ is uniquely arcwise connected. Let $g:S\rightarrow [0,\infty)$ be any function such that $g(s_0)=0$ and such that for every $s\in S$, the function $g$ is increasing on the arc $[s_0,s]$ of $S$. Then the function $d:S\times S\rightarrow [0,\infty)$, given by $d(s,t)=g(s)+g(t)-2g(s\wedge t)$, defines a metric on the set $S$.  Moreover, for every arc $[s,t]$ of the space $(S,\tau)$ the function $d(s,\cdot):[s,t]\rightarrow [0,d(s,t)]$ is an isometric embedding: $|d(s,u)-d(s,v)|=d(u,v)$ for all $u,v\in [s,t]$. (See \cite[pp.409--411]{MO} for details.)
\end{remark}

\begin{remark} If there is any arc $[s,t]$ in the space $(S,\tau)$ of Remark~\ref{ex} on which the given function $g$ is not continuous, then $[s,t]$ is not an arc of the metric space $(S,d)$. This logical pitfall appears to have been overlooked by the authors of \cite{MO}, when in \cite[Theorem 4.9]{MO} they erroneously claimed convexity of the resulting metric while still holding back the additional assumption of local arcwise connectedness. (See also the first paragraph on p.397 of \cite{MO}.)
\end{remark}

Recall that we defined $\rho((\omega_n)_n,(\xi_n)_n)=\Big\|\overline{(\omega_n)_n}\Big\|+\Big\|\overline{(\xi_n)_n}\Big\|-2\Big\|\overline{(\omega_n)_n}\Cap \overline{(\xi_n)_n}\Big\|$.

\begin{theorem}\label{tree}
 The function \[d(\tilde{y},\tilde{z})=\rho\left(\overleftarrow{\varphi(\tilde{y}}),\overleftarrow{\varphi(\tilde{z})}\right)
 =\Big\|\overline{\overleftarrow{\varphi(\tilde{y})}}\Big\|+\Big\|\overline{\overleftarrow{\varphi(\tilde{z})}}\Big\|
 -2\Big\|\overline{\overleftarrow{\varphi(\tilde{y})}}\Cap\overline{\overleftarrow{\varphi(\tilde{z})}}\Big\|\] defines an $\mathbb{R}$-tree metric on $\tilde{X}$ which induces the given topology.
\end{theorem}

\begin{proof}

Based on  Remarks~\ref{increasing} and \ref{ex} and Corollary~\ref{meet}, the function  $d$ defines a metric on $\tilde{X}$.  It suffices to show that this metric induces the given topology on $\tilde{X}$, as this implies that $(\tilde{X},d)$ is an $\mathbb{R}$-tree by Remark~\ref{ex}. We use Corollaries~\ref{coarc} and~\ref{meet} to describe the metric $d$:
 Let $\tilde{y},\tilde{z}\in \tilde{X}$ and choose any parametrizations $\tilde{\beta}:[0,1]\rightarrow \tilde{X}$ of $[\tilde{x},\tilde{y}]$,  $\tilde{\gamma}:[0,1]\rightarrow \tilde{X}$  of $[\tilde{x},\tilde{z}]$, and  $\tilde{\alpha}:[0,1]\rightarrow \tilde{X}$ of $[\tilde{x},\tilde{y}\wedge\tilde{z}]$.
 Put  $\beta=q\circ \tilde{\beta}=(\beta_n)_n$, $\gamma=q\circ \tilde{\gamma}=(\gamma_n)_n$ and $\alpha=q\circ \tilde{\alpha}=(\alpha_n)_n$.
 Then $\tilde{y}=[\beta]$,  $\tilde{z}=[\gamma]$, $\tilde{y}\wedge \tilde{z}=[\alpha]$ and \begin{eqnarray*} d(\tilde{y},\tilde{z})&=&\|(\omega_n(\beta_n))_n\|+\|(\omega_n(\gamma_n))_n\|-2\|(\omega_n(\alpha_n))_n\|\\
 &=&\lim_{n\rightarrow \infty}\Big( |\omega_n(\beta_n)|+|\omega_n(\gamma_n)|-2|\omega_n(\alpha_n)|\Big)
 \end{eqnarray*}
Put $y=(y_n)_n=q\circ\tilde{y}=\beta(1)$ and let $f_{n,\infty}:X\rightarrow X_n$ denote coordinate projection.

  First, let $[\delta]\in \tilde{X}$, let $U$ be an open subset of $X$ and suppose that $\tilde{y}\in B([\delta],U)$.  We wish to find an $\epsilon>0$ such that if $d(\tilde{y},\tilde{z})<\epsilon$, then $\tilde{z}\in B([\delta],U)$. Notice that $B([\delta],U)=B([\beta],U)$. Choose $n$ sufficiently large so that the combinatorial 6\nobreakdash-neighborhood $E_n$ of $y_n$ in $X_n$ is such that  $y\in f^{-1}_{n,\infty}(E_n)\subseteq U$.
   In order to prove that $\tilde{z}\in B([\beta],U)$, it suffices to show that $q_n([\tilde{z},\tilde{y}])\subseteq E_n$.
 To this end, write $\omega_n(\beta_n)=v_1v_2\cdots v_k/\ast$. Put $\epsilon=\min\{|v_1|/2^n,|v_2|/2^n,\dots,|v_k|/2^n\}$. Suppose that $d(\tilde{y},\tilde{z})<\epsilon$. We wish to show that $q_n([\tilde{z},\tilde{y}])=q_n([\tilde{z},\tilde{y}\wedge\tilde{z}])\cup q_n([\tilde{y}\wedge\tilde{z},\tilde{y}])\subseteq E_n$. Let $F_n$ be the combinatorial 3-neighborhood of $y_n$ in $X_n$. Then   we must have $q_n([\tilde{y}\wedge\tilde{z},\tilde{y}])\subseteq F_n\subseteq E_n$. (Otherwise,  $d(\tilde{y},\tilde{z})\geqslant \|(\omega_n(\beta_n))_n\|-\|(\omega_n(\alpha_n))_n\|\geqslant |v_{k-2}|\geqslant 2^n\epsilon>\epsilon$, by Lemma~\ref{three}; a contradiction.) Let $H_n$ be the combinatorial 3-neighborhood of $q_n(\tilde{y}\wedge \tilde{z})$ in $X_n$. Then  $q_n([\tilde{z},\tilde{y}\wedge\tilde{z}])\subseteq H_n\subseteq E_n$. (Otherwise, we have $\omega_n(\alpha_n)=v_1v_2\cdots v_s/\ast$ and $\omega_n(\gamma_n)=v_1v_2\cdots v_su_1u_2\cdots u_m/\ast$ with $s\leqslant k$ and $m\geqslant 3$. Then, by Lemmas~\ref{three} and~\ref{one}, $d(\tilde{y},\tilde{z})\geqslant \|(\omega_n(\gamma_n))_n\|-\|(\omega_n(\alpha_n))_n\|\geqslant\left(|v_1|+|v_2|+\cdots + |v_s|+|u_1|\right)-\left(|v_1|+|v_2|+\cdots+|v_s|\right)= |u_1|\geqslant |v_s|/2^n\geqslant \epsilon$; a~contradiction.)

Next, let $\epsilon>0$ be given. We wish to find an open set $U\subseteq X$ such that $B([\beta],U)\subseteq
\{[\delta]\in\tilde{X}\mid d([\beta],[\delta])<\epsilon\}$. By Lemmas~\ref{one} and~\ref{two}, we may choose $n$  so that the weight of every letter of the $n^\text{th}$ word of every word sequence is less than $\epsilon/4$ and so that the length of every word sequence is within $\epsilon/8$ of the length of its $n^\text{th}$ word.
Let $E_n$ be any open vertex-star of $X_n$ with $y_n\in E_n$. Put $U=f^{-1}_{n,\infty}(E_n)$ and suppose $\tilde{z}\in B([\beta],U)$. Then $[\gamma]=[\beta\cdot \tau]$ with $\tau=(\tau_n)_n:[0,1]\rightarrow U$.
Let $\tilde{\tau}:([0,1],0)\rightarrow (\tilde{X},\tilde{y})$ be the lift of $\tau$ with $q\circ \tilde{\tau}=\tau$. Then $\tilde{\tau}(1)=[\beta\cdot\tau]=\tilde{z}$. Since $\tilde{X}$ is uniquely arcwise connected, we have $[\tilde{y},\tilde{z}]\subseteq \tilde{\tau}([0,1])$. Since $\tilde{X}$ is simply connected (Lemma~\ref{CurtisFort} and Remark~\ref{simply}),  we may assume that $\tau$ equals the image of a parametrization of the arc $[\tilde{y},\tilde{z}]$ under the mapping $q$. Since the image of $\tau_n$ lies in  $E_n$ and since $E_n$ contains only one vertex of $X_n$, we get see that $\omega_n(\beta_n)$ and $\omega_n(\alpha_n)$ differ by at most one letter and so do $\omega_n(\gamma_n)$ and $\omega_n(\alpha_n)$. Since the weight of such a letter is less than $\epsilon/4$, we get from our choice of $n$ that $d(\tilde{y},\tilde{z})<\epsilon/8+\epsilon/8+2\epsilon/8+\epsilon/4+\epsilon/4 =\epsilon$.
\end{proof}

\begin{corollary}[Theorem~B]\label{pseudo} The function $\rho$, given by \[\rho((\omega_n)_n,(\xi_n)_n)=\Big\|\overline{(\omega_n)_n}\Big\|+\Big\|\overline{(\xi_n)_n}\Big\|-2\Big\|\overline{(\omega_n)_n}\Cap \overline{(\xi_n)_n}\Big\|,\]
 defines a pseudo metric on $\overleftarrow{\mathcal R}$ such that $\rho((\omega_n)_n,(\xi_n)_n)=0 \Leftrightarrow (\omega_n)_n\doteq (\xi_n)_n$.\linebreak  The resulting metric space $(\dot{\overleftarrow{\mathcal R}},\rho)$ is an $\mathbb{R}$-tree.
\end{corollary}

\begin{proof}
It remains to show that $\rho((\omega_n)_n,(\xi_n)_n)=0 \Leftrightarrow (\omega_n)_n\doteq (\xi_n)_n$.

If $(\omega_n)_n\doteq (\xi_n)_n$, then $\overline{(\omega_n)_n}=\overline{(\xi_n)_n}$ by Lemma~\ref{period}(c), so that we have $\rho((\omega_n)_n,(\xi_n)_n)=0$. Conversely, say $\rho((\omega_n)_n,(\xi_n)_n)=0$ with $(\omega_n)_n,(\xi_n)_n\in \overleftarrow{\mathcal R}$. Let $\tilde{y},\tilde{z}\in \tilde{X}$ be the unique elements with $(\omega_n)_n\doteq \overleftarrow{\varphi(\tilde{y})}$ and $(\xi_n)_n\doteq \overleftarrow{\varphi(\tilde{z})}$. Then $d(\tilde{y},\tilde{z})=0$ by Lemma~\ref{period}(c). Hence,  $\tilde{y}=\tilde{z}$ by Theorem~\ref{tree}, so that $(\omega_n)_n\doteq (\xi_n)_n$.
\end{proof}

\begin{theorem}[Theorem~E]\label{lc}
 If for every $n$, the essential multiplicity of every letter $v\in V_n$ is finite, then the quotient space $\dot{\overleftarrow{\mathcal R}}/\overleftarrow{\mathcal G}$ is homeomorphic to $X$. The essential multiplicity of every letter is finite if and only if $X$ is locally path connected.
\end{theorem}

\begin{proof}
First suppose that the essential multiplicity of every letter is finite. We will show that $X$ is locally path-connected, which implies that  $\dot{\overleftarrow{\mathcal R}}/\overleftarrow{\mathcal G}$ is homeomorphic to $X$ \cite[Theorem~4.10(c)]{FZ2}. It suffices to show that $X$ is locally connected. Following \cite{GM}, we show that $X$ has Property S: {\em every open cover of $X$ can be refined by a finite cover of connected subsets of $X$}. Fix $n$ and let $\{C_1, C_2, \dots, C_m\}$ be the collection of all closed vertex stars of $X_n$. Let $U_i$ be the open combinatorial 1\nobreakdash-neighborhood of $C_i$ in $X_n$. Since every vertex of $X_n$ has finite  essential multiplicity, there is a $K>n$ such that for every $k\geqslant K$ the number of components of $(f_n\circ f_{n+1}\circ\cdots\circ f_{k+1})^{-1}(U_i)$ which intersect $(f_n\circ f_{n+1}\circ\cdots\circ f_{k+1})^{-1}(C_i)$ is constant; say, this number is $l$. Label these components $V_{k,i}^1, V_{k,i}^2, \dots, V_{k,i}^l$ so that $f_k(V_{k+1,i}^j)\subseteq V_{k,i}^j$ for all $i=1,2,\dots, m$ and $j=1, 2, \dots, l$. Since both $i$ and $j$ range over finite index sets and since the sets $V_{k,i}^j$ cover $X_k$, we see that the sets
$S_i^j=$\[\lim_{\longleftarrow}\left(X_1 \stackrel{f_1}{\longleftarrow}  \cdots \stackrel{f_{K-2}}{\longleftarrow} X_{K-1} \stackrel{f_{K-1}}{\longleftarrow}  cl(V_{K,i}^j)\stackrel{f_K}{\longleftarrow} cl(V_{K+1,i}^j) \stackrel{f_{K+1}}{\longleftarrow} cl(V_{K+2,i}^j) \stackrel{f_{K+2}}{\longleftarrow} \cdots\right)\]  cover $X$. Also, each $S_i^j$ is connected.
  Note that given any open cover $\mathcal U$ of $X$, we can choose $n$ sufficiently large so that $\{S_i^j\mid i=1,2,\dots, m; j=1, 2, \dots, l\}$ refines~$\mathcal U$.

Now suppose $X$ is locally connected and let  $v\in V_n$.
Let $Star(v,X_n)$ denote the open star of the vertex $v$ in $X_n$. Then $c_k(v)$ equals the number of components of $(f_n\circ f_{n+1}\circ \cdots\circ f_{k})^{-1}(Star(v,X_n))$ which intersect $(f_n\circ f_{n+1}\circ \cdots\circ f_{k})^{-1}(\{v\})$.
  Let $f_{k+1,\infty}:X\rightarrow X_{k+1}$ denote coordinate projection, i.e., $f_{k+1,\infty}((x_i)_i)=x_{k+1}$.\linebreak
    By \cite[Lemma~2]{GM}, the numbers $c_k(v)$ are bounded by the (finite) number of (open) components of $(f_n\circ f_{n+1}\circ \cdots\circ f_{k}\circ f_{k+1,\infty})^{-1}(Star(v,X_n))$ which intersect the (compact) set $(f_n\circ f_{n+1}\circ \cdots\circ f_k\circ f_{k+1,\infty})^{-1}(\{v\})$. Hence, the essential multiplicity of $v$ is finite.
\end{proof}

\section*{Acknowledgements}
The authors gratefully acknowledge the support of this research by the Faculty Internal Grants Program of Ball State University and the Polish Ministry of Science (KBN 0524/H03/2006/31).

\end{document}